\documentclass{amsart}
\usepackage{amssymb, amsmath, amsfonts, amscd}

\input xypic

\theoremstyle{plain}
\newtheorem{theorem}{Theorem}[section]
\newtheorem{corollary}[theorem]{Corollary}
\newtheorem{lemma}[theorem]{Lemma}
\newtheorem{proposition}[theorem]{Proposition}

\newtheorem{example}[theorem]{Example}

\theoremstyle{definition}
\newtheorem{definition}[theorem]{Definition}

\theoremstyle{remark}

\numberwithin{equation}{theorem}

\newcommand{\J}{\mathcal{J}}

\renewcommand{\O}{\mathcal{O} }

\newcommand{\Hom}{\operatorname{Hom}}

\newcommand{\Der}{\operatorname{Der} }
\newcommand{\End}{\operatorname{End} }
\newcommand{\Spec}{\operatorname{Spec} }
\renewcommand{\P}{\operatorname{P} }
\renewcommand{\H}{\operatorname{H} }
\newcommand{\D}{\operatorname{D}^1_f(A) }
\newcommand{\Dl}{\operatorname{D} }
\newcommand{\iD}{\iota}

\newcommand{\U}{\operatorname{U}} 
\newcommand{\Uo}{\operatorname{U}^{\otimes} }

\newcommand{\Diff}{\operatorname{Diff}}

\newcommand{\Z}{\operatorname{Z} }
\newcommand{\tL}{\tilde{L} }
\newcommand{\tLone}{\tilde{L_1}}

\newcommand{\ta}{\tilde{\alpha} }
\newcommand{\tp}{\tilde{\pi} }

\newcommand{\LR}{\underline{\text{LR}(A/k)}}

\newcommand{\Conn}{\underline{\operatorname{Conn}} }
\newcommand{\Mod}{\underline{\operatorname{Mod} }}
\newcommand{\Ext}{\operatorname{Ext}}
\newcommand{\ext}{\operatorname{ext} }

\newcommand{\Pic}{\operatorname{Pic} }

\newcommand{\p}{\partial_1}
\newcommand{\pp}{\partial_2}

\newcommand{\ppp}{\partial_3}

\newcommand{\C}{\operatorname{C} }

\newcommand{\fal}{ \alpha^*(f) }

\begin{document}

\title{Extensions of Lie algebras of differential operators}

\author{Helge \"{O}ystein Maakestad}

\email{\text{h\_maakestad@hotmail.com}}

\thanks{}


\date{January 2019}

\begin{abstract} In this paper we introduce the notion  $\Dl$-Lie algebra and prove some elementary properties of $\Dl$-Lie algebras, the category of $\Dl$-Lie algebras, the category of connections on a 
$\Dl$-Lie algebra and extensions of $\Dl$-Lie algebras.  A $\Dl$-Lie algebra is a Lie-Rinehart algebra equipped with an $A\otimes_k A$-module structure
and a canonical central element $\iD$ and a compatibility property between the $k$-Lie algebra structure and the $A\otimes_k A$-module structure. 
Several authors have studied non-abelian extensions of Lie algebras, super Lie algebras, Lie algebroids and holomorphic Lie algebroids
and we give in this note an explicit constructions of all non-abelian extensions a $\Dl$-Lie algebra $\tL$ by an $A$-Lie algebra $(W,[,])$ where $\tL$ is projective as left $A$-module and $W$ is an $A\otimes_k A$-module
with $IW=0$ for $I$ the kernel of the multiplication map. As a corollary we get an explicit construction of all non-abelian extensions of an Lie-Rinehart algebra $(L,\alpha)$ by an $A$-Lie algebra $(W,[,])$ where $L$ is projective as left $A$-module. 
\end{abstract}

\maketitle

\tableofcontents

\section{Introduction}   The aim of this paper is to introduce a new structure - a D-Lie algebra - and to study non-abelian extensions of D-Lie algebras, non-abelian extensions of Lie-Rinehart algebras and a new 
characteristic class for  vector bundles. In the first section of the paper I introduce non-abelian extensions of Lie-Rinehart algebras and I use this construction to define a new characteristic class $nc_1(E)$ of an arbitrary
$A$-module $E$ with a connection $\nabla: L \rightarrow \End_k(E)$ where $L$ is a Lie-Rinehart algebra. The class $nc_1(E)$ is independent of choice of connection $\nabla$, and this construction gives a new characteristic class
for any finite rank vector bundle on any affine algebraic variety, including stably trivial vector bundles. If $E$ is a stably trivial vector bundle on an affine algebraic variety and  if $c_i(E)$ is a characteristic class satisfying the Whitney sum formula, it follows $c_i(E)=0$ for all $i\geq 1$. Hence for a stably trivial vector bundle $E$ it follows the Chern class $c_i(E)=0$ is zero. Hence if we consider the Chern class of $E$ with values in groups such as the 
Chow group, the Grothendieck group, the deRham cohomology group or the  Lie-Rinehart cohomology group, it follows Chern classes of stably trivial vector bundles are zero. As an example: if $X$ is a regular affine hypersurface and if $E$ is  the  the tangent bundle or cotangent bundle of $X$, it follows $E$  is stably trivial, hence all Chern classes $c_i(E)=0$ are zero.  The new class $nc_1(E)$ introduced in section 1 is non trivial for such vector bundles: I prove that $nc_1(T_{S^2})\neq 0$ is non trivial for the tangent bundle $T_{S^2}$ on the real 2-sphere $S^2$. The classical Chern classes $c_i(T_{S^2})=0$ are zero for the tangent bundle on the real 2-sphere, hence the new characteristic class $nc_1(E)$ is stronger than the classical Chern class. If $nc_1(E)=0$ is trivial it follows all Chern classes are trivial.

Let $A$ be an arbitrary commutative unital $k$-algebra where $k$ is an arbitrary commutative untal ring. The \emph{module of first order differential operators} of $A/k$ is defined as the set
of $k$-linear endomorphisms $\partial : A\rightarrow A$ such that for any $a\in A$ it follows the operator $\partial \circ \phi_a-\phi_a \circ \partial$ is multiplication with an element in $A$. Here $\phi_a$ is the 
operator on $A$ given by multiplication with the element $a$. We denote
the set of such operators by $\Dl^1_0(A)$. The ring $\End_k(A)$ is canonically a left $A\otimes_k A$-module and
$k$-Lie algebra  and it follows $\Dl^1_0(A)$ is an $A\otimes_k A$-submodule and $k$-Lie subalgebra of $\End_k(A)$.  There is moreover an isomorphism $\Dl^1_0(A)\cong A\oplus \Der_k(A)$ as $A\otimes_k A$-modules.
If $f\in \Z^2(\Der_k(A),A)$ is a 2-cocycle we may define the following: Let $\D:=A\oplus \Der_k(A)$
with the following $k$-Lie product: 
\[  [u,v]:=(x(b)-y(a)+f(x,y), [x,y]) \]
where $u:=(a,x),v:=(b,y)\in \D$. Define the following left and right $A$-module structure on $\D$: Let 
\[ cu:=c(a,x):=(ca,cx),\text{ and }uc:=(a,x)c:=(ac+x(c),cx) \]
for $c\in A$. It follows $\D$ is a left $A\otimes_k A$-module and $k$-Lie algebra. There is a canonical map $\pi: \D \rightarrow \Der_k(A)$ defined by $\pi(a,x):=x$ and $\pi$ is a map of left $A\otimes_k A$-modules
and $k$-Lie algebras. There is a central element $z:=(1,0)\in \D$ with $\pi(z)=0$. The following holds for any $u,v\in \D$ and $c\in A$:
\begin{align}
&\label{d1}  [u,cv]=c[u,v]+\pi(u)(c)v 
\end{align}
and
\begin{align}
&\label{d2} uc=cu+\pi(u)(c)z .
\end{align}
Hence $(\D, \pi)$ with its structure as left $A$-module is an Lie-Rinehart algebra. A $\Dl$-Lie algebra is a generalization of $\D$. We say a 5-tuple $(\tL, \ta,\tp,[,], \iD)$ is a \emph{$\Dl$-Lie algebra}
if the following holds: $\tL$ is a left $A\otimes_k A$-module and $k$-Lie algebra. The map $\tp: \tL \rightarrow \D$ is a map of $A\otimes_k A$-modules and $k$-Lie algebras where $f\in Z^2(\Der_k(A),A)$ 
is a 2-cocycle. The map $\tp: \tL \rightarrow \Der_k(A)$ is a map of $A\otimes_k A$-modules and $k$-Lie algebras where we have given $\Der_k(A)$ the trivial right $A$-module structure. 
The element $\iD\in Z(\tL)$ is a central element, and the following holds for any $u,v\in \tL$ and $c\in A$:
\begin{align}
&\label{dl1} [u,cv]=c[u,v]+\tp(u)(c)v 
\end{align}
and 
\begin{align}
&\label{dl2} uc=cu+\tp(u)(c)\iD .
\end{align}
Hence the if we view the pair $(\tL, \tp)$ as a left $A$-module and a map of left $A$-modules and $k$-Lie algebras, it follows $(\tL, \tp)$ is an Lie-Rinehart algebra. 
By Lemma \ref{modulelie} the following gives the relationship between the Lie product and the right $A$-module structure on $\tL$:
\begin{align}
&\label{dl3} [u,vc]=c[u,v]+\tp(u)(c)v+\tp(u) \circ \tp(v)(c) \iD.
\end{align}

The main motivation for the introduction of a $\Dl$-Lie algebra is the following. Given an arbitrary $L$-connection $(E,\nabla)$ with the property that $\nabla$ is of \emph{curvature type $f$} for a 2-cocycle $f\in \Z^2(L,A)$,
there is a generalized universal enveloping algebra  $\U(A,L,f)$ (see \cite{maa1}). The associative ring $\U(A,L,f)$ has the property that there exists an exact equivalence $\psi_f$ between the  category of left $\U(A,L,f)$-modules and maps of modules and the  category  of $L$-connections of curvature type $f$ and morphisms of $L$-connections. The equivalence $\psi_f$ preserves injective and projective objects.
In the paper \cite{maa1} I introduced Ext and Tor groups of connections of curvature type $f$ using the algebra $\U(A,L,f)$ and the equivalence $\psi_f$. A connection of curvature type $f$ is a special case of a non-flat connection. In the paper \cite{maa145} I generalize this construction further  to the case of an arbitrary $L$-connection $(E,\nabla)$ using the notion $\Dl$-Lie algebra. I define the Ext and Tor groups
of any pair of connections $(E,\nabla), (F,\nabla')$ on any Lie-Rinehart algebra using the notion $\Dl$-Lie algebra $\tL$ and the \emph{universal ring} $\Uo(\tL)$ of a $\Dl$-Lie algebra. In \cite{maa145} I interpret
the category of connections  $\Conn(L)$ as the module category $\Mod(\Uo(\tL))$ of a $\Dl$-Lie algebra $\tL$ associated to $L$ and use such an interpretation to construct Ext and Tor groups. Since the Ext and Tor groups
are defined using the category of modules $\Mod(\Uo(\tL))$ on an associative ring $\Uo(\tL)$, it follows all functorial properties of these groups are immediate from the classical book of Cartan and Eilenberg on homological algebra (see \cite{cartan}). It would not have been possible to introduce Tor groups of arbitrary connections without the introduction of $\tL$, the universal ring $\Uo(\tL)$ and the interpretation of $\Conn(L)$ as the category of modules $\Mod(\Uo(\tL))$. One cannot use the Freyd-Mitchell Full Embedding Theorem (see \cite{freyd}) since the equivalence constructed in the Freyd-Mitchell Theorem does not preserve injective and projective objects in general. Hence the notion of a $\Dl$-Lie algebra $\tL$ and the corresponding universal  ring $\Uo(\tL)$ is essential if we want to study the Ext and Tor groups of arbitrary $L$-connections.

 In this paper we study connections on $\Dl$-Lie algebras and non-abelian extensions of $\Dl$-Lie algebras and Lie-Rinehart algebras. In the litterature one may 
find a treatment of the notion of non-abelian extensions of Lie algebras, Lie algebroids and holomorphic Lie algebroids (see  \cite{michor},\cite{giraud},  \cite{mackenzie} and \cite{tortella1}). 
Hence a holomophic version of the notion of a non-abelian extension have appeared in \cite{tortella1}. Non-abelian extensions of $\Dl$-Lie algebras have not appeared in the litterature since the notion $\Dl$-Lie algebra appears to be a new notion. I have not found it studied elsewhere.

In the main Theorem of the paper (Theorem \ref{maintheorem}) I do the following. Let $(\tL, \ta, \tp, [,],\iD)$ be a $\Dl$-Lie algebra where $\tL$ is projective as left $A$-module, and let $(W,[,])$ be an $A$-Lie algebra with 
$aw=wa$ for all $a\in A$ and $w\in W$. Let $S(\tL, (W,[,]))$ denote the set of all triples 
\begin{align}
&\label{triple} s:=(\nabla_s, \psi_s, \tilde{w})
\end{align}

where $\nabla_s:\tL\rightarrow \Der_k(W)$ is a connection and 
$\psi_s: \wedge^2 \tL \rightarrow W$ is a 2-cocycle satisfying the following property (property P): 
\begin{align}
&\label{propertyp} R_{\nabla_s}(x,y)=[\psi_s(x,y),-]\text{ and  }d^2_{\nabla_s}(\psi_s)=0.
\end{align}
Assume $\tilde{w}\in W$ is an element with 
\begin{align}
&\label{centralelement} \psi_s(\iD,-)=d^0_{\nabla_s}(\tilde{w})\text{ and }\nabla_s(\iD)(-)=ad(-\tilde{w}). 
\end{align}
Let two triples $s:=(\nabla_s, \psi_s, w_s)$ and $t:=(\nabla_t, \psi_t, w_t)$ be equivalent (written $s\equiv t$) if there is a map $b\in \Hom_A(\tL,W)$ such that
\begin{align}
 &\label{equiv1}\nabla_t(x)=\nabla_s(x)+[b(x),-] \\
 &\label{equiv2}\psi_t(x,y)=\psi_s(x,y)+d^1_{\nabla_s}(b)(x,y)+[b(x),b(y)]  \\
  &\label{equiv3} w_t=w_s-b(\iD) .
\end{align}
In Lemma \ref{equivS} we prove that the relation $\equiv$ defined by equations \ref{equiv1}, \ref{equiv2} and \ref{equiv3} is an equivalence relation on the set $S(\tL, (W,[,]))$.
Let $\ext^1(\tL, (W,[,])):= S(\tL, (W,[,]))/ \equiv$ be the quotient set. Let $\Ext^1(\tL, (W,[,]))  $ denote the set of equivalence classes of extensions of the $\Dl$-Lie algebra $\tL$ with the $A$-Lie algebra $(W,[,])$ 
modulo isomorphism of extensions. The main theorem (see Theorem \ref{maintheorem}) is the following:
\begin{theorem} Let $(\tL, \ta,\tp,[,], \iD)$ be a $\Dl$-Lie algebra with $\tL$ a projective $A$-module. Let $(W,[,])$ be an $A$-Lie algebra with $aw=wa$ for all $a\in A, w\in W$.
There is a one-to-one correspondence of sets
\[ F: \Ext^1(\tL, (W,[,])) \cong \ext^1(\tL, (W,[,])) \]
\end{theorem}

Hence the set $\ext^1(\tL, (W,[,]))$ of equivalence classes of triples $(\nabla, \psi,w)$ satisfying \ref{propertyp} and \ref{centralelement} classifies extensions of $\Dl$-Lie algebras.

In Theorem \ref{extconn} we prove the following result on connections on a classical Lie-Rinehart algebra $(L,\alpha)$: Let $(E,\nabla)$ be an $L$-connection. There is an Lie-Rinehart algebra
$(\End(L,(E,\nabla)), p_E)$ and  an extension $nc_1(E)$
 \begin{align}
&\label{nonab} 0 \rightarrow \End_A(E) \rightarrow \End(L,(E,\nabla)) \rightarrow L \rightarrow 0 
\end{align}
of Lie-Rinehart algebras - the \emph{non-abelian extension of $(L,\alpha)$ by  $(E,\nabla)$}.

\begin{theorem} Let $(L,\alpha)$ be an Lie-Rinehart algebra and let $\nabla: L\rightarrow \End_k(E)$ be a connection. 
Let $P\in \Hom_A(L,\End_A(E))$ and define $\nabla':=\nabla + P$. It follows $\nabla'$ is a connection on $E$ and there is a canonical isomorphism 
of non-abelian extensions $\End(L,(E, \nabla))\cong \End(L,(E, \nabla'))$ of Lie-Rinehart algebras.   Hence the extension class \ref{nonab} defined by  $\End(L,(E, \nabla))$ is 
independent of choice of connection $\nabla$. The $A$-module $E$ has a flat connection if and only if there is a map of of Lie-Rinehart algebras $s_E:L\rightarrow \End(L,(E,\nabla))$ with $p_E \circ s_E=Id_L$.
\end{theorem}

Hence the module $E$ has a flat connection if and only if the exact sequence
\[ 0 \rightarrow \End_A(E) \rightarrow \End(L,(E,\nabla)) \rightarrow L \rightarrow 0 \]
splits in the category of Lie-Rinehart algebras.

In Example \ref{real2sphere} we prove that the extension class $nc_1(T_{S^2})$ from \ref{nonab} is non trivial for the tangent bundle $T_{S^2}$ on the real 2-sphere $S^2$. The Chern classes $c_i(T_{S^2})=0$ are trivial, hence the class $nc_1(T_{S^2})$  is stronger than the Chern classes. It gives a non trivial characteristic class for any vector bundle on any algebraic variety, in particular for any stably trivial vector bundle on any algebraic variety. If $nc_1(E)=0$ it follows all Chern classes $c_i(E)=0$ are trivial.

Assume the ring $k$ contains a field of characteristic zero. 
Given an Lie-Rinehart algebra $(L,\alpha)$  and any connection $\nabla:L\rightarrow \End_k(E)$ where $E$ is a finite rank projective $A$-module we construct in Theorem \ref{mainext}
a non-trivial extension of $L(f)$ by $\End_A^{tr}(E)$. 

\begin{theorem} There is an extension of pre-$\Dl$-Lie algebras
\[0\rightarrow \End_A^{tr}(E) \rightarrow \End(L,(E,\nabla)) \rightarrow L(f) \rightarrow 0 .\]
Here  $f:=tr(R_\rho)$ is the trace of the curvature of $\rho$. The $A$-Lie algebra $\End_A^{tr}(E)$ is the quotient of $\End_A(E)$ by the sub-$A$-module
and $A$-Lie algebra of endomorphisms with trace zero. 
\end{theorem}

Hence non-abelian extensions of pre-$\Dl$-Lie algebras and $\Dl$-Lie algebras arise naturally when studying classical connections on finite rank projective modules.

We also give an explicit proof of the classification of all non-abelian extensions of an Lie-Rinehart algebra $(L,\alpha)$ by an $A$-Lie algebra $W$ where $L$ is a projective left $A$-module.
This classification has appeared in the paper \cite{tortella1} for holomorphic Lie algebroids and here I give an explicit proof of this classification in the algebraic case. Some of the proofs have already appeared
in the paper \cite{tortella1} hence there is nothing new about this classification. I include all details of the proof in this paper since I will need these details in the proof of the classification for extensions of $\Dl$-Lie algebras. 


\section{A new characteristic class for vector bundles}

In this section we give a complete classifcation of all non-abelian extensions of an $A$-Lie algebra $(W,[,])$ with an Lie-Rinehart algebra $(L,\alpha)$ where $L$ is a projective $A$-module.
This classification was previously known for $k$-Lie algebras, Lie-super algebras, Lie algebroids for differentiable manifolds and holomorphic Lie algebroids for complex manifolds 
(see \cite{michor}, \cite{mackenzie} and  \cite{tortella1}). I have not found algebraic proofs of these results published in the litterature and I need the results in the final chapter of the paper. This is why 
I include the results and proofs in complete detail. 

The aim of the section is to introduce non-abelian extensions of Lie-Rinehart algebras and use this construction to introduce a new characteristic class for vector bundles.
For any Lie-Rinehart algebra $(L, \alpha)$ and any connection $(E, \nabla)$, I construct a non-abelian extension class $nc_1(E)$ that is independent of choice of connection $\nabla$. The class is trivial 
if and only if $E$ has a flat connection. In \ref{real2sphere} I prove this class is non-trivial for the tangent bundle on the real two sphere. The Chern class of the tangent bundle on the real 2-sphere are trivial,
hence the class $nc_1$ is stronger than the Chern class.

Given two  Lie-Rinehart algebras $(L,\alpha), (L',\alpha')$ and a surjective map of Lie-Rinehart algebras $p:L' \rightarrow L$ is follows
the kernel $W:=ker(p)$ in a natural way is an $A$-Lie algebra. The $k$-linear bracket $[,]$ on $L'$ induce an $A$-linear bracket on $W$. If $L$ is a projective $A$-module there is a left  $A$-linear  section $s$
of $p$ and one may define the following maps: 
\begin{align}
&\label{nablas} \nabla_s:L\rightarrow \Der_k(W)\text{ by }\nabla_s(x)(w):=[s(x),w]
\end{align}
and 
\begin{align}
&\label{psis}\psi_s: \wedge^2 L\rightarrow W\text{ by }\psi_s(x,y):=[s(x),s(y)]-s([x,y]).
\end{align}
 It follows $\nabla_s$ is a connection on $W$ with curvature $R_{\nabla_s}(x,y)(w)=[\psi_s(x,y),w]$. Moreover $d^2_{\nabla_s}(\psi_s)=0$ where $d^2_{\nabla_s}$ is the Chevalley-Eilenberg differential for the connection 
$(W,\nabla_s)$.
Hence we may ask if it is possible to classify the set of extensions of the Lie-Rinehart algebra $(L,\alpha)$ by the $A$-Lie algebra $(W,[,])$ in terms of the maps $\nabla_s$ and $\psi_s$ and this is possible. 
The aim of this section is to do this for an arbitrary Lie-Rinehart algebra $(L,\alpha)$ where $L$ is an arbitrary projective $A$-module and $(W,[,])$ is any $A$-Lie algebra. 
In Corollary \ref{mainone} we classify such extensions  via pairs $(\nabla, \psi)$ where $\nabla$ is a connection, $\psi$ a 2-cocycle with $d^2_{\nabla}=0$ and $R_{\nabla}(x,y)=[\psi(x,y), -]$ modulo an equivalence relation on the set of such pairs $(\nabla,\psi)$ (see Definition \ref{equivLR}).  The classification is similar to the classification one gets for $k$-Lie algebras.

A Lie-Rinehart algebra is an algebraic version
of the notion of a Lie algebroid and much of the cohomological and homological constructions that can be proved for a $k$-Lie algebra over a field can be proved for an Lie-Rinehart algebra $L$
where $L$ is a projective $A$-module. The calculations are straightforward and can be done as a student exercise in homological algebra and extensions. Still I have not found it written out in detail in the litterature but it is probably well known to experts. The aim of this section is to do this explicitly and in the next section to generalize the classification to the case of non-abelian extensions of $\Dl$-Lie algebras.

In the following for a Lie-Rinehart algebra $(L,\alpha)$ and a connection $(W,\nabla)$, let 
\begin{align}
&\label{LRcomplexone} \C^i(L,(W,\nabla)):=\Hom_A(\wedge^i L,W) 
\end{align}
be the $i$'th \emph{Chevalley-Eilenberg complex} of the connection $(W,\nabla)$ with differential $d^i$. When $\nabla$ is a flat connection, the complex in \ref{LRcomplexone} is  referred to as the \emph{Lie-Rinehart complex}. This is due to Rineharts paper \cite{rinehart}, which was the first systematic study of the complex \ref{LRcomplexone}, the universal enveloping algebra $\U(A,L)$, the PBW theorem for $\U(A,L)$ and the cohomology group $\H^i(L,(W,\nabla))$ and homology group $\H_i(L,(W,\nabla))$ in the case when $L$ is a projective $A$-module and $(W,\nabla)$ a flat connection. For this reason many authors refer to the cohomology $\H^i(L,(W,\nabla))$ (resp. $\H_i(L,(W,\nabla))$) as the \emph{Lie-Rinehart cohomology} (resp. \emph{Lie-Rinehart homology}) of the flat connection $(W,\nabla)$. 

Many authors studied Lie-Rinehart algebras prior to Rineharts paper \cite{rinehart} (his PhD thesis) hence
the notion was well known to experts at the time. I do not have enough competence on the history of the development of the notion to write about this 
and interested readers should consult the book of Mackenzie \cite{mackenzie} and other papers by the same author. Interested readers shuld also consult the paper by Huebschmann \cite{huebschmann} which gives much information and references in section one.

The differential $d^i$ in the complex \ref{LRcomplexone} depends on the connection $\nabla$ and sometimes we will 
write $d^i_{\nabla}$ when it is unclear which connection on $W$ we speak about. Assume $[,]$ is an $A$-Lie product on the $A$-module $W$. 
See \cite{maa1} for some basic results on the complex \ref{LRcomplexone}.

\begin{proposition} \label{LRalgebra} Let $k$ be a commutative ring and $A$ a $k$-algebra. Let $(L;\alpha)$ and $(L',\alpha')$ be Lie-Rinehart algebras and let 
$p:L'\rightarrow L$ be a surjective map of Lie-Rinehart algebras. Assume $L$ is projective as left $A$-module and let $W:=Ker(p)$. It follows $(W,[,])$ is an $A$-Lie algebra.
Let $s: L\rightarrow L'$ be a left $A$-linear splitting of $p$ and let  $x\in L, w\in W$. Define $\nabla(x)(w):= [s(x),w]\in W$. It follows we get a connection
\[ \nabla:L\rightarrow \Der_k(W).\]
Let for $x,y\in L$ $\psi(x,y):=[s(x),s(y)]-s([x,y]).$ It follows we get an $A$-linear map
\[ \psi: \wedge^2 L\rightarrow W \]
with $d^2_{\nabla}(\psi)=0$ where $d^2_{\nabla}$ is the $2$'nd differential in the Lie-Rinehart complex of $(W,\nabla)$.
For $x,y\in L, w\in W$ it follows $R_{\nabla}(x,y)(w)=[\psi(x,y),w]$ hence $\nabla$ is a non-flat connection in general.
\end{proposition} 
\begin{proof}
The Proposition is a straight forward calculation.
\end{proof}

\begin{definition} \label{conditionA} Given an Lie-Rinehart algebra $(L,\alpha)$ and an $A$-Lie algebra $(W,[,])$ consider the following set $S((L,\alpha) ,(W,[,]))$ of pairs $(\nabla, \psi)$ where
$\nabla: L\rightarrow \Der_k(W)$ is a connection and $\psi \in \Z^2(L,W)$ is a 2-cocycle satisfying the following condition: Let $x,y \in L$ and $w\in W$:
\begin{align}
&\label{2set}R_{\nabla}(x,y)(w)=[\psi(x,y),w].
\end{align}
Let the condition in \ref{2set} be called \emph{condition A}.
\end{definition}

\begin{lemma} \label{change}Let $(L,\alpha)$ be an Lie-Rinehart algebra and let $(\nabla, \psi)$ be a pair satisfying condition A from Definition \ref{conditionA}. Let $P\in \Hom_A(L, \End_A(W))$ be an element with 
\[ P(x)([w,v])=[P(x)(w),v]+[w,P(x)(v)] \]
for all $w,v \in W$. It follows $\nabla':=\nabla+P$ is another connection 
\[ \nabla' : L\rightarrow \Der_k(W).\]
There is for all $x,y\in L, w\in W$ an equality
\[ R_{\nabla'}(x,y)(w)=[\psi(x,y),w]+d^2_{ad\nabla}(P)(x,y)(w)+[P(x),P(y)](w) .\]
\end{lemma}
\begin{proof} The proof is a straightforward calculation.
\end{proof}

Hence by Lemma \ref{change} it follows that given a pair $(\nabla, \psi)$ satisfying condition A, when we add a potential $P$ to $\nabla$ to get another connection $\nabla':=\nabla + P$ we may
not always find a 2-cocycle $\psi'\in \Z^2(L,W)$ such that the new pair $(\nabla', \psi')$ satisfies condition A.

Let $(L,\alpha)$ be an Lie-Rinehart algebra and let $(W,[,])$ be an $A$-Lie algebra. Let $(\nabla, \psi)\in S((L,\alpha), (W,[,]))$ be an element.
Define the following bracket on the direct sum $W\oplus L$:
\begin{align}
\label{bracket}  [(w,x),(v,y)]:= ([w,v]+\nabla(x)(v)-\nabla(y)(w)+\psi(x,y), [x,y]).
\end{align}
for $(w,x),(v,y)\in W\oplus L$.
Define the map
\[ \alpha_W : W\oplus L \rightarrow \Der_k(A) \]
by $\alpha_W(w,x):=\alpha(x)\in \Der_k(A)$.
Let $(L(W,(\nabla,\psi)), \alpha_W)$ denote the left $A$-modue $W\oplus L$ equipped with the bracket $[,]$ from \ref{bracket}, and the map $\alpha_W$.
There is a canonical map 
\[ p: L(W,(\nabla,\psi)) \rightarrow L \]
defined by
\[ p(w,x):=x\in L.\]
We get an exact sequence of left $A$-modules
\begin{align}
&\label{LRW} 0 \rightarrow W \rightarrow L(W,(\nabla,\psi)) \rightarrow L \rightarrow 0.
\end{align}

\begin{lemma}\label{LRalgebra} Let $(W,[,])$ be an $A$-Lie algebra and let $(L,\alpha)$ be an Lie-Rinehart algebra. Let $(\nabla,\psi) \in S((L,\alpha),(W,[,]))$ be an element.
It follows the pair 
\[ (L(W,(\nabla, \psi)), \alpha_W) \]
is an Lie-Rinehart algebra. The exact sequence \ref{LRW} is an extension of Lie-Rinehart algebras.
\end{lemma}
\begin{proof} Let us prove the Jacobi identity for the bracket defined in \ref{bracket}. The rest is an exercise for the reader.
Let $(u,x),(v,y),(w,z)\in W\oplus L$. Consider the sum
\begin{align}
\label{sum} [(u,x),[(v,y),(w,z)]]+[(v,y),[(w,z),(u,x)]]+[(w,z),[(u,x),(v,y)]] .
\end{align}

For the first term in \ref{sum} we get
\[ [(u,x),[(v,y),(w,z)]]=\]
\[ [u[vw]]+[u,\nabla(y)(w)]-[u,\nabla(z)(v)]+[u,\psi(y,z)]+\nabla(x)([v,w])+\nabla(x)\nabla(y)(w) \]
\[ -\nabla(x)\nabla(z)(w) +\nabla(x)(\psi(y,z))-\nabla([y,z])(u)+\psi(x,[y,z]), [x,[y,z]]) .\]
For the second term in the sum \ref{sum} we get
\[ [(v,y),[(w,z), (u,x)]]=\]
\[ ([v,[w,u]]+[v,\nabla(z)(u)]-[v,\nabla(x)(w)]+[v,\psi(z,x)]+\nabla(y)([w,u])+\nabla(y)\nabla(z)(u) \]
\[ -\nabla(y)\nabla(x)(w)+\nabla(y)(\psi(z,x))-\nabla([z,x])(v)+\psi(y,[z,x]), [y,[z,x]]) \]
For the third term in the sum \ref{sum} we get
\[ [(w,z),[(u,x),(v,y)]]=\]
\[ ([w,[u,v]]+[w,\nabla(x)(v)]-[w,\nabla(y)(u)]+[w,\psi(x,y)]+\nabla(z)([u,v])+\nabla(z)\nabla(x)(v) \]
\[ -\nabla(z)\nabla(y)(u)+\nabla(z)(\psi(x,y)) -\nabla([x,y])(w)+\psi(x,[y,z]),[z,[x,y]]).\]
We get
\[ [u,\psi(y,z)]+[v,\psi(z,x)]+[w,\psi(x,y)]+\]
\[ R_{\nabla}(x,y)(w)+R_{\nabla}(z,x)(v)+R_{\nabla}(y,z)(u)\]
\[ -[\psi(y,z),u]-[\psi(z,x),v]-[\psi(x,y), w] =0.\]
We get
\[ \nabla(x)(\psi(y,z))-\nabla(y)(\psi(z,x))+\nabla(z)(\psi(x,y)) \]
\[ -\psi([x,y],z) +\psi([x,z],y)-\psi([y,z],x)=d^2_{\nabla}(\psi)(x,y,z)=0.\]
We get
\[ [(u,x),[(v,y),(w,z)]]+[(v,y),[(w,z),(u,x)]]+[(w,z),[(u,x),(v,y)]]=0 \]
and the Jacobi identity holds. The rest is an exercise for the reader and the Lemma follows.
\end{proof}


\begin{definition} \label{equivLR}
Define the following relation on the set $S((L,\alpha), (W,[,]))$: Given two elements $(\nabla,\psi),(\nabla',\psi') \in S((L;\alpha), (W,[,]))$. We say $(\nabla,\psi) \equiv (\nabla',\psi')$
if and only if there is a map $b\in \Hom_A(L,W)$ such that 
\[ \nabla'(x)=\nabla(x)+[b(x), -] \]
and 
\[ \psi'(x,y)=\psi(x,y)+d^1_{\nabla}(b)(x,y)+[b(x),b(y)] \]
for all $x,y\in L$.
\end{definition}

\begin{lemma} \label{relation} The relation $\equiv$ on the set $S((L,\alpha), (W,[,]))$ is an equivalence relation.
\end{lemma}
\begin{proof} 
Clearly $(\nabla,\psi)\equiv(\nabla,\psi)$ with $b=0$. Assume there is an element $b\in \Hom_A(L,W)$ with 
\[ \nabla'(x)=\nabla(x)+[b(x),-] \]
and
\[ \psi'(x,y)=\psi(x,y)+d^1_{\nabla}(x,y)+[b(x),b(y)].\]
It follows
\[ \nabla(x)=\nabla'(x)+[(-b)(x),-]:=\nabla'(x)+[a(x),-] \]
where $a:=-b\in \Hom_A(L,W)$.
We want to prove that
\[ \psi(x,y)=\psi'(x,y)+d^1_{\nabla'}(a)(x,y)+[a(x),a(y)].\]
A calculation shows that
\[ d^1_{\nabla'}(a)(x,y)=-d^1_{\nabla}(b)(x,y)-2[b(x),b(y)].\]
We get by assumption
\[ \psi'(x,y)+d^1_{\nabla'}(a)(x,y)+[a(x),a(y)]=\]
\[\psi(x,y)+d^1_{\nabla}(b)(x,y)+[b(x),b(y)]-d^1_{\nabla}(b)(x,y)-2[b(x),b(y)]+[b(x),b(y)]=\psi(x,y).\]
Hence
\[ \psi(x,y)=\psi'(x,y)+d^1_{\nabla'}(a)(x,y)+[a(x),a(y)] \]
and it follows
\[ (\nabla',\psi')\equiv (\nabla, \psi)\]
hence the relation $\equiv$ is symmetric.
If there are elements $a,b\in \Hom_A(L,W)$ with 
\[ \nabla'(x)=\nabla(x)+[b(x),-] \]
and
\[ \psi'(x,y)=\psi(x,y)+d^1_{\nabla}(x,y)+[b(x),b(y)]\]
and
\[ \nabla''(x)=\nabla'(x)+[a(x),-]  \]
and
\[ \psi''(x,y)=\psi(x,y)'+d^1_{\nabla}(a)(x,y)+[a(x),a(y)]\]
it follows
\[ \nabla''(x)=\nabla(x)+[(a+b)(x),-] \]
and
\[ \psi''(x,y)=\psi(x,y)+d^1_{\nabla}(a+b)(x,y)+[(a+b)(x),(a+b)(y)] \]
and $a+b\in\Hom_A(L,W)$. The Lemma follows.
\end{proof}

\begin{definition} Let $(L,\alpha)$ be an Lie-Rinehart algebra and let $(W,[,])$ be an $A$-Lie algebra. Let $\ext^1((L,\alpha), (W,[,])):= S((L,\alpha), (W,[,]))/\equiv$
be the quotient set.
\end{definition}

Note: There is in general no structure of abelian group on the set 
\[ \ext^1((L,\alpha), (W,[,])). \]

Let 
\begin{align}
\label{exact}0  \rightarrow W \rightarrow L' \rightarrow L \rightarrow 0 
\end{align}

be an exact sequence of Lie-Rinehart algebras with $p:L'\rightarrow L$  a map of Lie-Rinehart algebras and $W:=ker(p)$.  It follows $(W,[,])$ is an $A$-Lie algebra in a canonical way. 
Let  

\[ T(L,W):=T((L,\alpha),(W,[,])) \]

be the set of short exact sequences of Lie-Rinehart algebras as in \ref{exact} where $(W,[,])$ is a fixed $A$-Lie algebra and $(L,\alpha)$ a fixed Lie-Rinehart algebra.
Let two extensions 
\[ L' \rightarrow^{p'} L\text{ and }L''\rightarrow^{p''} L \]
 in $T(L,W)$ be equivalent if there is an isomorphism of Lie-Rinehart algebras $\rho:L' \rightarrow L''$ making all diagrams commute.
We write $(L',p')\equiv (L'',p'')$ if $(L',p')$ and $(L'', p'')$ are equivalent extensions. The relation $\equiv$ is an equivalence relation on $T(L,W)$ since the map $\rho$ is an isomorphism. 

\begin{definition}
Let 
\[ \Ext^1((L,\alpha), (W,[,]) ):= \Ext^1(L,W):=T(L,W)/\equiv \]
 be the quotient set of $T(L,W)$ with respect to the equivalence relation $\equiv$.
\end{definition}

Assume $L$ is a projective left  $A$-module and $s:L\rightarrow L'$ a left splitting of $p$. Let $t$ be another splitting of $p$. It follows $t=s+b$ with $b\in \Hom_A(L,W)$. Let $(\nabla_s,\psi_s), (\nabla_t,\psi_t)$ be the 
two elements in $S((L,\alpha), (W,[,]))$ defined by  $s$ and $t$.  
Define for the splitting $s$ the  exact sequence of left $A$-modules as done in Proposition \ref{LRalgebra}:
\begin{align}
&\label{LRext} 0 \rightarrow W \rightarrow L(W,(\nabla_s,\psi_s)) \rightarrow L \rightarrow 0 .
\end{align}
It follows from Proposition \ref{LRalgebra} the sequence \ref{LRext} is an extension of Lie-Rinehart algebras.

\begin{lemma} \label{lemmaiso} For any splitting $s$ of $p$ it follows the pair $(L(W,(\nabla_s,\psi_s)), \alpha_W)$ is an Lie-Rinehart algebra. The sequence
\[ 0 \rightarrow W \rightarrow L(W,(\nabla_s,\psi_s)) \rightarrow L \rightarrow 0\]
is an exact sequence of Lie-Rinehart algebras.
If $t$ is another left $A$-linear splitting of $p$ with $s=t+b$ with $b\in \Hom_A(L,W)$ it follows there is an isomorphism
\[ \rho: L(W, (\nabla_s,\psi_s)) \rightarrow L(W, (\nabla_t,\psi_t) ) \]
of $A$-modules and $k$-Lie algebras defined by
\[ \rho(w,x)=(w+b(x),x) \]
making all diagrams commute. Hence the map $\rho$ is an isomorphism of extensions of Lie-Rinehart algebras.
\end{lemma}
\begin{proof} The first two claims follow from Proposition \ref{LRalgebra}.

Let $(w,x),(v,y)\in W\oplus L$. We get
\[ \rho([(w,x),(v,y)])=([w,v]+\nabla_s(x)(v)-\nabla_s(y)(w)+\psi_s(x,y)+b([x,y]),[x,y]).\]
We get
\[ [\rho(w,x),\rho(v,y)]=\]
\[ ([w,v]+\nabla_t(x)v)+[(b(x),v]-(\nabla_t(y)(w)+[b(y),w])+\nabla_t(x)(b(y))- \]
\[ \nabla_t(y)(b(x))+\psi_t(x,y)+[b(x),b(y)],[x,y]).\]
It follows
\[ [\rho(w,x),\rho(v,y)]-\rho([(w,x),(v,y)])=\]
\[ (\nabla_t(x)(v)+[b(x),v]-\nabla_s(x)(v)-(\nabla_t(y)(w)+[b(y),w]-\nabla_s(y)(w)) + \]
\[ d^1_{\nabla_t}(b)(x,y)+\psi_t(x,y)-\psi_s(x,y),0).\]
Since $\nabla_s(x)=\nabla_t(x)+[b(x),-]$ and $\psi_s(x,y)=\psi_t(x,y)+d^1_{\nabla_t}(b)(x,y)+[b(x),b(y)]$ we get
\[ \rho([(w,x),(v,y)])=[\rho(w,x),\rho(v,y)] \]
hence $\rho$ is an isomorphism of Lie-Rinehart algebras. It follows $\rho$ is an isomorphism of extensions and the Lemma follows.
\end{proof}

\begin{definition} \label{defs}
Let $(L,\alpha)$ be a Lie-Rinehart algebra and let $(W,[,])$ be an $A$-Lie algebra. Let $s:=(\nabla,\psi)$ be a a pair where $\nabla:L\rightarrow \Der_k(A)$ is a connection and 
$\psi: \wedge^2 L \rightarrow W$ a 2-cocycle satisfying condition A. 
Let $L(W,s):=L(W,(\nabla,\psi))$  and let  $(L(W,s ), [,], \alpha_W)$ be the triple from Lemma \ref{LRalgebra}. 
Let
\[ p_s: L(W,s) \rightarrow L \]
be the map $p_s(w,x):=x\in L$.
\end{definition}
By definition  $(L(W,s), [,], \alpha_W)$  is an Lie-Rinehart algebra and the map $p_s$ is a surjective map of Lie-Rinehart algebras.

\begin{theorem} \label{mainthm}  Let $s:=(\nabla,\psi)$ and $s':=(\nabla',\psi')$ be two pairs satisfying condition A and consider the Lie-Rinehart algebras $L(W,s)$ and $L(W, s')$ from 
Definition \ref{defs}.
Define the map $\rho: L(W,s) \rightarrow L(W,s')$ by $\rho(w,x)=(w+b(x),x)$ where $b\in \Hom_A(L,W)$. It follows $\rho$ is an isomorphism of left $A$-modules.
Consider the conditions
\begin{align}
&\label{iso1}\nabla(x)(v)=\nabla'(x)(v)+[b(x),v] \\
&\label{iso2}\psi(x,y)=\psi'(x,y)+d^1_{\nabla'}(b)(x,y)+[b(x),b(y)]
\end{align}
It follows the map $\rho$ is a map of $k$-Lie algebras if and only if conditions \ref{iso1} and \ref{iso2} hold.
\end{theorem}
\begin{proof} Let $u:=(w,x), v:=(v,y)\in L(W,s):=W\oplus L$. We get
\[ \rho([u,v])= \]
\[([w,v]+\nabla(x)(v)-\nabla(y)(w)+\psi(x,y)+b([x,y]),[x,y]).\]
We moreover get
\[ [\rho(u),\rho(v)]=\]
\[([w,v]+\nabla'(x)(v)+[b(x),v]-(\nabla'(y)(w)+[b(y),w]) +\nabla'(x)(b(y)) \]
\[ -\nabla'(y)(b(x))+\psi'(x,y)+[b(x),b(y)],[x,y]).\]
It follows
\[  [\rho(u),\rho(v)]-\rho([u,v])=\]
\[ (\nabla'(x)(v)+[b(x),v]-\nabla(x)(v) -(\nabla'(y)(w)+[b(y),w]-\nabla(y)(w) ) \]
\[ +d^1_{\nabla'}(b)(x,y)+\psi'(x,y)-\psi(x,y)+[b(x),b(y)], 0) \]
Hence if $\rho$ is a map of Lie algebras for all $u,v$ it follows if we let $y=0$, we get
\[ \nabla'(x)(v)+[b(x),v]=\nabla(x)(v) \]
for all $x\in L$ and $v\in W$. It follows
\[ [\rho(u),\rho(v)]-\rho([u,v])= \]
\[ (\psi'(x,y)+d^1_{\nabla'}(b)(x,y)+[b(x),b(y)]-\psi(x,y),0).\]
Hence since $\rho$ is a map of Lie algebras it follows
\[ \psi(x,y)=\psi'(x,y)+d^1_{\nabla'}(b)(x,y)+[b(x),b(y)] \]
and the first claim follows. 
Assume equations \ref{iso1} and \ref{iso2} holds. By the above calculation it follows $\rho$ is a map of Lie algebras. The Theorem follows.
\end{proof}

\begin{lemma} Let $(L,\alpha)$ be an Lie-Rinehart algebra where $L$ is a projective $A$-module and let $(W,[,])$ be an $A$-Lie algebra. Let 
\[ 0\rightarrow W \rightarrow L' \rightarrow L \rightarrow 0 \]
be an exact sequence of Lie-Rinehart algebras with $p:L' \rightarrow L$ a surjective map of Lie-Rinehart algebras with a left $A$-linear section $s$. Let $\nabla:L\rightarrow \Der_k(W)$
be the connection associated to $s$ and $\psi: \wedge^2 L \rightarrow W$ the 2-cocycle associated to $s$. It follows $(\nabla,\psi)\in S((L,\alpha),(W,[,]))$.
Define the following map of sets:
\[ F: \Ext^1((L,\alpha), (W,[,]) ) \rightarrow \ext^1((L,\alpha), (W,[,])) \]
by
\[ F(p:L'\rightarrow L):= (\nabla , \psi) \operatorname{mod} (\equiv)  .\]
The map  $F$ is a well defined map of sets. Here $\equiv$ is the equivalence relation from Definition \ref{equivLR} defined on the set $S((L,\alpha),(W,[,]))$.
\end{lemma}
\begin{proof} We must prove that if $p_1:L_1\rightarrow L$ and $p_2:L_2\rightarrow L$ are two equivalent extensions then
\[ F(p_1:L_1\rightarrow L ) \equiv F( p_2:L_2\rightarrow L ) \]
in the quotient $\ext^1((L,\alpha), (W,[,]))$. Chose two $A$-linear sections $s_1$ of $p_1$ and $s_2$ of $p_2$ and define the connections and 2-cocycles $ (\nabla_1 ,\psi_1 )$ and
$( \nabla_2 , \psi_2 )$. From Lemma \ref{lemmaiso} we get isomorphisms of extensions
\[ L(W, s_1)  \cong L_1\text{ and }L(W, s_2)  \cong L_2.\]
If $\rho: L_1 \rightarrow L_2$ is an isomorphism of extensions we get an induced isomorphism of extensions
\[ \tilde{\rho}: L(W,s_1) \rightarrow L(W, s_2) .\]
It follows from Theorem \ref{mainthm}  that there is an element $b\in \Hom_A(L,W)$ such that
\[ \nabla_1(x)=\nabla_2(x)+[b(x), -] \]
and
\[ \psi_1(x,y)=\psi_2(x,y)+d^1_{\nabla_2}(x,y)+[b(x),b(y)] \]
It follows there is an equality 
\[  F( p_1:L_1\rightarrow L ) \equiv F( p_2:L_2\rightarrow L ) \]
in the quotient set $\ext^1((L,\alpha),(W,[,]))$,  hence the map $F$ is well defined and the Lemma is proved.
\end{proof}

\begin{corollary} \label{mainone} Let $(L,\alpha)$ be an Lie-Rinehart algebra where $L$ is projective as left $A$-module. Assume $(W,[,])$ is an $A$-Lie algebra. There is a one-to-one correspondence
of sets
\[ F: \Ext^1((L,\alpha), (W,[,]) ) \cong \ext^1((L,\alpha),(W,[,])).\]
Hence the set $\ext^1((L,\alpha),(W,[,]))$ classifies extensions of $(L,\alpha)$ by $(W,[,])$.
\end{corollary}
\begin{proof}
Assume we are given two extensions $p_1:L_1\rightarrow L $ and $p_2:L_2\rightarrow L$ with 
\[ F( p_1:L_1\rightarrow L ) \equiv F( p_2:L_2\rightarrow L ) \text{ in $\ext^1((L,\alpha), (W,[,]))$}. \]
If $ F( p_1:L_1\rightarrow L )=(\nabla_1,\psi_1)$ and  $F(p_2:L_2\rightarrow L ) =( \nabla_2,\psi_2 )$ it follows from Theorem \ref{mainthm} there is an element $b\in \Hom_A(L,W)$
with 
\[ \nabla_2(x) = \nabla_1(x) +[b(x), -] \]
and
\[ \psi_2(x,y)=\psi_1(x,y)+d^1_{\nabla_1 }(b)(x,y)+[b(x),b(y)].\]
It follows from Theorem \ref{mainthm} the map
\[ \tilde{\rho}: L(W, (\nabla_1,\psi_1) )\rightarrow L(W, (\nabla_2,\psi_2)) \]
defined by
\[ \tilde{\rho}(w,x):=(w+b(x),x) \]
is an isomorphism of extensions. It follows $L_1\cong L_2$ is an isomorphism of extensions, hence the map $F$ is an injection. 

Assume $(\nabla, \psi)\in \ext^1((L,\alpha),(W,[,]))$ is an element.
Define the following Lie product on $W\oplus L$ as done in Proposition \ref{LRalgebra}:
\[ [(w,x),(v,y)]:=([w,v]+\nabla(x)(v)-\nabla(y)(w)+\psi(x,y), [x,y]).\]
It follows from Proposition \ref{LRalgebra} the $A$-module $L(W, (\nabla,\psi))$ is an Lie-Rinehart algebra and furthermore $L(W, (\nabla,\psi))$ is an extension of $L$ by $W$. 
Since $L$ is a projective $A$-module it follows there is a section $s$ of the natural projection morphism
$p:L(W, (\nabla,\psi))\rightarrow L$ and $s(x)=(b(x),x)$ where $b\in \Hom_A(L,W)$. One checks that
$F(L(W, (\nabla,\psi), p)=( \nabla_s , \psi_s)$ where
\[ \nabla_s(x)=\nabla(x)+[b(x), -] \]
and
\[ \psi_s(x,y)=\psi(x,y)+d^1_{\nabla}(b)(x,y)+[b(x),b(y)] .\]
Hence
\[ F(L(W, (\nabla,\psi), p):=(\nabla_s,\psi_s )\equiv (\nabla, \psi)  \]
and the map $F$ is surjective. The Corollary is proved.
\end{proof}

Note: When $W$ is an abelian Lie algebra we recover the classical classification of abelian extensions of $L$ by $W$ as discussed in \cite{maa1}.

\begin{example} The non-abelian extension associated to a connection.\end{example}

In the classical situation of a connection on a Lie-Rinehart algebra we get in a canonical way a non-abelian extension associated to a connection. Let $(L,\alpha)$ be an Lie-Rinehart algebra.
If $\nabla: L\rightarrow \End_k(E)$ is an arbitrary  connection we get an induced connection
\[ \tilde{\nabla}: L\rightarrow \Der_k(\End_A(E)) \]
defined by
\[ \tilde{\nabla}(x)(\phi):=[\nabla(x),\phi].\]
One checks that for any pair $\phi,\psi\in \End_A(E)$ it follows 
\[ \tilde{\nabla}(x)([\phi,\psi])=[\tilde{\nabla}(x)(\phi),\psi]+[\phi, \tilde{\nabla}(x)(\psi)].\]
One gets for any elements $x,y \in L$ the following formula:
\[R_{\tilde{\nabla}}(x,y)(\phi)=[R_{\nabla}(x,y),\phi] \]
and $R_{\nabla}\in \Z^2(L,\End_A(E))$ hence we get since $\End_A(E)$ is an $A$-Lie algebra a non abelian extension
\begin{align}
&\label{extE} 0 \rightarrow \End_A(E) \rightarrow \End(L,(E,\nabla)) \rightarrow L \rightarrow 0 
\end{align}
where $p_E:\End(L,(E,\nabla))\rightarrow L$ is the canonical map.
Here $\End(L,(E,\nabla)):=\End_A(E) \oplus L$ has the following Lie product: Given $(\phi,x),(\psi,y)\in \End(L,E)$ we define

\[ [(\phi,x),(\psi,y)]:= ([\phi,\psi]+\tilde{\nabla}(x)(\psi)-\tilde{\nabla}(y)(\phi)+R_{\nabla}(x,y),[x,y]).\]
One checks $\End(L,(E,\nabla))$ is an Lie-Rinehart algebra.
Define 
\[ \alpha_E: \End(L,(E,\nabla))\rightarrow \Der_k(A) \]
  by $\alpha_E(\phi,x):=\alpha(x)\in \Der_k(A)$.

\begin{theorem} \label{extconn} Let $(L,\alpha)$ be an Lie-Rinehart algebra and let $\nabla: L\rightarrow \End_k(E)$ be a connection. 
Let $P\in \Hom_A(L,\End_A(E))$ and define $\nabla':=\nabla + P$. It follows $\nabla'$ is a connection on $E$.
There is a canonical isomorphism of non-abelian extensions $\End(L,(E,\nabla))\cong \End(L,(E,\nabla'))$, hence the extension class \ref{extE} in $\Ext^1((L,\alpha), (\End_A(E),[,]))$ 
is independent of choice of connection $\nabla$. 
The following holds: The $A$-module $E$ has a flat connection if and only if there is a map of of Lie-Rinehart algebras $s_E:L\rightarrow \End(L,E,\nabla)$ with $p_E \circ s_E=Id_L$.
\end{theorem}
\begin{proof} Let $\nabla:L\rightarrow \End_k(E)$ and let $ad\nabla: L\rightarrow \Der_k(\End_A(E))$ be the adjoint connection defined by $ad\nabla(x)(\phi):=[\nabla(x),\phi]$. Let $\nabla':=\nabla+P$. It follows
\[ ad\nabla'(x)=ad\nabla(x)+[P(x),-] .\]
Moreover
\[ R_{\nabla'}(x,y)=R_{\nabla}(x,y)+d^1_{ad\nabla}(P)(x,y)+[P(x),P(y)] \]
hence by Theorem \ref{mainthm} it follows there is an isomorphism of extensions
\[ \End(L,(E,\nabla))\cong \End(L,(E,\nabla')) \]
and the first claim is true. We prove the second claim: A splitting $s_E:L\rightarrow \End(L,(E,\nabla)$) is on the following form: $s(x)=(\rho(x),x)$ where $\rho\in \Hom_A(L,\End_A(E))$. It follows
\[ [s(x),s(y)]=[(\rho(x),x),(\rho(y),y)]=\]
\[ ([\rho(x),\rho(y)]+[\nabla(x),\rho(y)]-[\nabla(y),\rho(x)]+R_{\nabla}(x,y),[x,y])=\]
\[ ([\rho(x),\rho(y)]+[\nabla(x),\rho(y)]+[\rho(x),\nabla(y)]+[\nabla(x),\nabla(y)]-\nabla([x,y]),[x,y]).\]
It follows $s([x,y])=[s(x),s(y)]$ if and only if the following holds:
\[ [\rho(x),\rho(y)]+[\nabla(x),\rho(y)]+[\rho(x),\nabla(y)] +[\nabla(x),\nabla(y)]-\nabla([x,y])-\rho([x,y])=0\]
which is if and only if
\[ [\nabla(x)+\rho(x),\nabla(y)+\rho(y)]-(\nabla+\rho)([x,y])=R_{\nabla+\rho}(x,y)=0.\]
Hence $p_E$ has a splitting which is $A$-linear and a map of $k$-Lie algebras if and only if there is a map $\rho \in \Hom_A(L,\End_A(E))$ such that $\nabla':=\nabla +\rho$ is a flat connection.
The Theorem follows.
\end{proof}

\begin{definition} Let $nc_1(E)$ be the cohomology class in $\Ext^1((L,\alpha), \End_A(E))$ defined by the sequence \ref{extE}.
\end{definition}

Note: The class $nc_(E)$ is defined for an arbitrary connection $(E,\nabla)$. By Theorem \ref{extconn} it follows the class $nc_1(E)$ is independent of choice of connection $\nabla$.
When $E$ is finite rank and projective, there is a trace map $tr: \End_A(E)\rightarrow A$ and when $L$ is projective as $A$-module, this map gives a map (see \cite{maa1})
\[ \phi: \Ext^1((L,\alpha), \End_A(E)) \rightarrow \Ext^1((L,\alpha),A)\cong \H^2(L,A) \]
and it follows $\phi(nc_1(E))=c_1(E)\in \H^2(L,A)$ where $c_1(E)$ is the well known first Chern class of $E$.

\begin{example} \label{examplenonab} Atiyah sequences for connections on Lie-Rinehart algebras \end{example}

Hence the $A$-module $E$ has a flat connection if and only if the exact sequence \ref{extE} is split in the category of Lie-Rinehart algebras. Hence we may use the set
$\Ext^1((L,\alpha), (\End_A(E),[,]) )$ and Corollary \ref{mainone} to study the problem of determining if an $A$-module $E$ has a flat connection. Recall from \cite{maa1}, Lemma A.19 there is for any Lie-Rinehart algebra
$(L,\alpha)$ and any $A$-module $E$ an exact sequence of $A\otimes_k A$-modules
\begin{align}
&\label{atiyahseq}0\rightarrow E  \rightarrow J^1_L(E) \rightarrow L\otimes_A E   \rightarrow 0 
\end{align}
with the property that an $A\otimes_k A$-linear section $s: L\otimes_A E \rightarrow J^1_L(E) $ corresponds to a connection $\nabla: L\rightarrow \End_k(E)$.
The sequences \ref{extE} and \ref{atiyahseq} are defined for an arbitrary Lie-Rinehart algebra $(L,\alpha)$ and an arbitrary $A$-module $E$.

\begin{example} \label{exnonab} The non-abelian extension by a flat connection.\end{example}

Let $k$ be a field of characteristic zero and let $A:=k[x,1/x]$, $L_A:=\Der_k(A)$.  Let $R:=k[x_i, 1/x_i]$ for $i=1,..,n$ and let $L_R:=\Der_k(R)$.

\begin{lemma} \label{derham} There is an isomorphism 
\begin{align}
&\label{drone} \H^1(L_R,R) \cong k\{\frac{1}{x_i }dx_i :\text{ with $1\leq i  \leq n$} \}\\
&\label{drtwo}\H^2(L_R,R) \cong k\{\frac{1}{x_ix_j }dx_i \wedge dx_j :\text{ with $1\leq i < j  \leq n$} \}
\end{align}
of vector spaces over $k$.
\end{lemma}
\begin{proof} The Lemma follows from the Kunneth formula and the fact that $\H^0(L_A,A)=k$ and $\H^1(L_A,A)=k\frac{1}{x}dx$.
\end{proof}



Let $E:=R^d:=R\{e_1,..,e_d\}$ be the free $R$-module of rank $d$ and let $\rho:L\rightarrow \End_k(E)$ be defined by $\rho(x)(\sum c_ie_i):=\sum_i x(c_i)e_i$.
It follows $\rho$ is a flat connection.  Let $\phi:=(\phi_{ij}), \psi:=(\psi_{ij})\in \End_A(E)$. It follows we get a non-abelian extension
\begin{align}
&\label{nonabflat} 0 \rightarrow \End_A(E) \rightarrow \End(L, (E,\rho)) \rightarrow L \rightarrow 0
\end{align}
with the following Lie product: Let $u:=(\psi,x),v:=(\psi,y)\in \End(L,(E,\rho))$ and define
\[ [u,v]:=([\phi, \psi]+[\rho(x),\psi]-[\phi, \rho(y)]+R_{\rho}(x\wedge y),[x,y]) =\]
\[ ([\phi, \psi]+(x(\psi_{ij}))- (y(\phi_{ij})) ,[x,y]),\]
since $\rho$ is a flat connection. It follows the sequence \ref{nonabflat} splits in the category of Lie-Rinehart algebras.

\begin{example} Non-abelian extensions and non-triviality of projective modules. \end{example}

In \cite{maa0} the following is proved: Let $k$ be any field and let $f(x_1,..,x_n) \in k[x_1,..,x_n]$ be a regular hypersurace in the sense that the Jacobi ideal
$J(f)=(1)$ is the unit ideal. Let $A:=k[x_i]/(f)$ and let $\Omega:=\Omega^1_{A/k}$ be the module of Kahler differentials on $A$. It follows $\Omega$ is a finite rank projective 
$A$-module and there is a projective basis $\{B,B^*\}$ for $\Omega$ and a corresponding algebraic connection

\[ \nabla_B: \Der_k(A) \rightarrow \End_A(\Omega) \]

induced by $\Omega$.  We may construct the Chern classes 

\[  c_i(\Omega), c_i(\Der_k(A)) \in \H^{2i}_{DR}(A) \]

for all $i \geq 1$ where $\H^i_{DR}(A)$ is algebraic deRham cohomology of $A$.  The following holds:

\begin{theorem}  There is for every $i\geq 1$ an equality $c_i(\Omega)= c_i(\Der_k(A))=0$.   \end{theorem}
\begin{proof}  See \cite{maa0} Section 4.
\end{proof}

Hence for any affine regular hypersurface $H:=Spec(A)$ it follows its tangent and cotangent bundle has all Chern classes in algebraic deRham cohomology equal to zero. Hence the Chern class in algebraic deRham chomology is a weak invariant in the sense that it does not detect that these vector bundles are non-trivial.

\begin{example} \label{real2sphere} Vector bundles on the real $n$-sphere and characteristic classes. \end{example}

If $k$ is the field of real numbers and if $S^n$ is the real n-sphere, it follows the real tangent bundle $T_{S^n}$ is known to be topologically non-trivial, but all higher Chern classes vanish since it is stably trivial.
Hence the Chern class in the Chow group, Grothendieck group, algebraic deRham cohomology and Lie-Rinehart cohomology does not detect the non-triviality of $T_{S^n}$. Any theory of characteristic classes satisfying the Whitney sum formula has to be trivial for $T_{S^2}$. The characteristic class $nc_1(-)$ is non trivial for such vector bundles, hence the class is stronger than such theories of characteristic classes.

Let $f:=x^p+y^q+z^r-1$ and let $A:=k[x,y,z]/(f)$. with $S:=\Spec(A)$. In the paper \cite{maa0} I construct explicit formulas for non-flat algebraic connections

\[ \nabla: \Der_k(A) \rightarrow \End_A(\Omega) \]

with trace of curvature $tr(K_{\nabla})=0$. The set of alll connections on $\Omega$ equals $\Hom_A(\Der_k(A), \End_A(\Omega))$ and this is a large set. It is not clear if $\Omega$ has a flat connection.  The formulas are explicit and look as follows: Let $\Der_k(A)$ be generated by the derivations $\partial_1, \partial_2, \partial_3$. We get the following
formulas:

Let $\nabla(\partial_1)$ be the following operator:
\[ \nabla(\p)=D_{\p}+\p(M) .\]
One checks that the following holds:
\begin{align}
&\label{PARTIALone} \p(M)=
\end{align}
\[
\begin{pmatrix} -pqx^{p-1}y^{q-1}  & -p(p-1)x^{p-2}y^q +\frac{p^2}{q}x^{2(p-1)}  & -\frac{qp(p-1)}{r}x^{p-2}y^{q-1}z   \\
                -\frac{q^2}{p}y^{2(q-1)}+q(q-1)x^py^{q-2} &  pqx^{p-1}y^{q-1} & \frac{pq(q-1)}{r}x^{p-1}y^{q-2}z \\
 -\frac{qr}{p}y^{q-1}z^{r-1}  & \frac{pr}{q}x^{p-1}z^{r-1}  &  0 
\end{pmatrix}.
\]

One checks that  $\pp(M)$ is the following matrix:
\begin{align}
&\label{PARTIALtwo} \pp(M)=
\end{align}

\[
\begin{pmatrix} -prx^{p-1}z^{r-1} & -\frac{rp(p-1)}{q}x^{p-2}yz^{r-1} & -p(p-1)x^{p-2}z^r+\frac{p^2}{r}x^{2(p-1)} \\
        -\frac{qr}{p}y^{q-1}z^{r-1} & 0 & \frac{pq}{r}x^{p-1}y^{q-1} \\
     -\frac{r^2}{p}z^{2(r-1)}+r(r-1)x^pz^{r-2} & \frac{pr(r-1)}{q}x^{p-1}yz^{r-2} & prx^{p-1}z^{r-1}
\end{pmatrix}.
\]

One checks that $\ppp(M)$ is the following matrix:
\begin{align}
&\label{PARTIALTHREE} \ppp(M)=
\end{align}
\[
\begin{pmatrix} 0 &  -\frac{pr}{q}x^{p-1}z^{r-1} & \frac{pq}{r}x^{p-1}y^{q-1} \\
       -\frac{rq(q-1)}{p}xy^{q-2}z^{r-1} & -qry^{q-1}z^{r-1} & -q(q-1)y^{q-2}z^r+\frac{q^2}{r}y^{2(q-1)} \\
    \frac{qr(r-1)}{p}xy^{q-1}z^{r-2} & -\frac{r^2}{q}z^{2(r-1)}+r(r-1)y^qz^{r-2} & qry^{q-1}z^{r-1}
\end{pmatrix}.
\]

It follows we get well defined first order differential operators
\[ D_{\partial_i}+\partial_i(M):\Omega \rightarrow \Omega,\]
for $i=1,2,3$ giving explicit formulas for an algebraic connection

\[ \nabla: \Der_k(A) \rightarrow \End_A(\Omega).\]

The connection $\nabla$ has by \cite{maa0} trace of curvature $tr(R_{\nabla})=0$ equal to zero, hence the corresponding Chern class $c_1(\Omega)$ is zero. If $k$ is the real numbers and if $p=q=r=2$ it follows $\Omega$ is the cotangent bundle of the real 2-sphere which is known to be topologically non-trivial since it is the dual ot the tangent bundle. 

There is by this section a non-abelian cohomology class

\[  nc_1(\Omega) \in \Ext^1(\Der_k(A), \End_A(\Omega)) \]

which is independent with respect to choice of connection $\nabla$. This class is trivial if and only if $\Omega$ has a flat connection. If $\Omega$ is trivial it has a flat connection hence if the class $nc_1(\Omega)$ is non-trivial we can conclude that $\Omega$ is non-trivial. Hence the non-abelian cohomology class $nc_1(\Omega)$ can be used to decide if the module $\Omega$ is non-trivial. 

In \cite{milnor}, page 294 the following is mentioned without proof: If $M$ is a simply connected smooth manifold and $E$ is a smooth vector bundle with a flat connection, it follows $E$ is a trivial as smooth vector bundle. A proof of this may be found in \cite{taubes}.
Any algebraic connection on $T_{S^2}$ is smooth, hence if $T_{S^2}$ has a flat algebraic connection $\nabla$, it follows $\nabla$ is is a flat and smooth connection on $T_{S^2}$. This implies that $T_{S^2}$ is topologically trivial which is a contradiction. Hence $T_{S^2}$ does not have a flat algebraic connection. From this it follows the characteristic class $nc_1(T_{S^2})$  is non-trivial for $T_{S^2}$.
Hence the class introduced above detects that the tangent bundle on the real 2-sphere is non trivial. It follows the class introduced above is stronger than the Chern class since the Chern class is trivial for $T_{S^2}$.




\section{Non-abelian extensions of $\Dl$-Lie algebras}

In this section we introduce the notion of a $\Dl$-Lie algebra, the category of $\Dl$-Lie algebras and the category of modules on $\Dl$-Lie algebras. The construction is modeled on the module of first order
$k$-linear differential operators $\Dl^1_0(A)$ and the notion of an Lie-Rinehart algebra.

Let $k$ be a commutative ring and let $A$ be a commutative $k$-algebra. Let $\Dl^1_0(A)$ be the module of first order $k$-linear differential opreators on $A$.
There is an isomorphism $i:\Dl^1_0(A)\cong A\oplus \Der_k(A)$ given by $i(\partial):=(\partial(1), \partial-\partial(1))$ and there is a canonical inclusion map
\[ j: A\oplus \Der_k(A) \rightarrow \End_k(A) \]
defined by
\[ j(a,\delta):= \phi_a + \delta \]
where $\phi_a(b):=ab$ is multiplication with $a$. It follows $\Dl^1_0(A)$ has a canonical $A\otimes_k A$-module structure and $k$-Lie algebra structure 
induced by the inclusion $\Dl^1_0(A) \subseteq \End_k(A)$. It is  defined as follows:
\[ a(b,\delta):=(ab,a\delta)\]
\[(b,\delta)a:=(ba+\delta(a),a\delta) \]
and
\[ [(b,\delta),(c,\eta)] := (\delta(c)-\eta(b),[\delta,\eta]) .\]

Note that there is no non-trivial right $A$-module structure on $\Der_k(A)$. We must consider the abelian extension $\Dl^1_0(A)$ to get something non-trivial since $\delta \circ \phi_a = \phi_{\delta(a)} +a\delta$ in
$\End_k(A)$.

Given a 2-cocycle $f\in \Z^2(\Der_k(A),A)$ we may construct a new $k$-Lie algebra structure on $A\oplus \Der_k(A)$ as follows:
\begin{align}
&\label{lie1} [(b,\delta), (c,\eta)]:= (\delta(c)-\eta(b)+f(\delta,\eta), [\delta,\eta]) .
\end{align}
Let $(\D ,[,])$ denote the $A\otimes_k A$ module $A\oplus \Der_k(A)$ equipped with the $k$-Lie bracket defined in \ref{lie1}. Let $\pi: \D \rightarrow \Der_k(A)$ be the canonical projection map.
It follows $\pi$ is $A\otimes_k A$-linear and a map of $k$-Lie algebras. Let $\overline{z}:=(1,0)\in \D$. It follows $[\overline{z},u]=0$ for all elements $u\in \D$ hence the element 
$\overline{z}\in Z(\D)$ is in the center $Z(\D)$ of the Lie algebra $\D$.

\begin{lemma}  \label{modulelie} Let $u:=(a,\delta), v:=(b,\eta)\in \D$ and let $c\in A$. The bracket $[,]$ from \ref{lie1} satisfies the following:
\begin{align}
&\label{d0}  vc=cv+\pi(v)(c)\overline{z}
\end{align}
\begin{align}
&\label{d1}  [u,cv]=c[u,v]+\pi(u)(c)v 
\end{align}
and
\begin{align}
&\label{c2} [u,vc]=[u,cv]+\pi(u)\pi(v)(c) \overline{z}. 
\end{align}
\end{lemma}
\begin{proof} Claim \ref{d0} and \ref{d1} follows from an explicit calculation. Let us verify Claim \ref{c2}: 
We get
\[ vc:=(b,\eta)c:=(bc+\eta(c),c\eta)=c(b,\eta)+\eta(c)\overline{z}=cv+\pi(v)(c)\overline{z}.\]
It follows
\[ [u,vc]=[u,cv]+[u, \pi(v)(c)\overline{z}]= \pi(v)(c)[u,\overline{z}]+\pi(u)\pi(v)(c)\overline{z} \]
hence
\[ [u,vc]=[u,cv]+\pi(u)\pi(v)(c)\overline{z}.\]
\end{proof}

Note: Condition \ref{d1} is the well known anchor-condition from the definition of a classical Lie-Rinehart algebra.

The notion of a $\Dl$-Lie algebra is an extension of the module of first order differential operators $\Dl^1_0(A)$ or more generally the module $\D$:

\begin {definition} \label{dLie} Let $A$ be a commutative $k$-algebra where $k$ is a commutative ring.  A $5$-tuple $(\tL, \ta, \tp, [,],\iD)$ is a \emph{$\Dl$-Lie algebra} if the following holds:
$(\tL,[,])$ is a $k$-Lie algebra and $A\otimes_k A$-module.  The map $\ta:\tL\rightarrow \D$ is a map of $k$-Lie algebras and $A\otimes_k A$-modules where $f\in \Z^2(\Der_k(A),A)$ is a 2-cocycle. 
The map $\tp: \tL\rightarrow \Der_k(A)$ is a map of $k$-Lie algebras and $A\otimes_k A$-modules with $\pi \circ \ta = \tp$. We give $\Der_k(A)$ the trivial right $A$-module structure.
Let $u,v\in \tL$ and $c\in A$. The element $\iD\in Z(\tL)$ is a central element with $\tp(\iD)=0$. 
The following holds:
\begin{align}
&\label{dl1} uc = cu+\tp(u)(c) \iD 
\end{align} 
\begin{align}
&\label{dl2} [u,cv]=c[u,v]+\tp(u)(c) v 
\end{align}

Given two $\Dl$-Lie algebras $(\tL, \ta, \tp, [,], \iD)$ and $(\tilde{L'}, \tilde{\alpha'},\tilde{\pi'},[,], \iD')$. A \emph{map of $\Dl$-Lie algebras} is a map
\[ \phi: \tL \rightarrow \tilde{L'} \]
of $A\otimes_k A$-modules and $k$-Lie algebras with the property that $\phi(\iD)=\iD' $ and with $\tilde{\alpha'}\circ \phi =\ta$. 
A $\Dl$-Lie algebra  $\tL$ is also referred to as a \emph{Lie algebra of first order differential operators on $A/k$}.
Let $\underline{\Dl-\text{Lie}}$ denote the category of $\Dl$-Lie algebras and maps of $\Dl$-Lie algebras.

A \emph{pre-$\Dl$-Lie algebra} is a 4-tuple $(\tL, \tp, [,], \iD)$ where $(\tL, [,])$ is a $k$-Lie algebra and $A\otimes_k A$-module, the map $\tp: \tL \rightarrow \Der_k(A)$ is a map
of $k$-Lie algebras and $A\otimes_k A$-modules. The element $\iD\in Z(\tL)$ is a central element with $\tp(\iD)=0$. The Lie product and central element $\iD$ satisfies equations \ref{dl1} and \ref{dl2}.
\end{definition}

Note: the 5-tuple $(\D, id, \pi, [,], z)$ where $z:=(1,0)$ is a $\Dl$-Lie algebra. Given a $\Dl$-Lie algebra  $(\tL, \ta, \tp, [,], \iD)$ it follows $(\tL, \tp)$ is an ordinary Lie-Rinehart algebra. Hence a $\Dl$-Lie algebra 
is an Lie-Rinehart algebra equipped with the structure of a right $A$-module with a certain compatibility between the Lie bracket and the right $A$-module structure given in Lemma \ref{rightcomp}.
There is moreover the canonical central element $\iD$ and the relationship between $\iD$ and the $A\otimes_k A$-module structure on $\tL$.

\begin{example} Variants of Lie algebras in the litterature. \end{example}

Let $k$ be any commutative ring and let $A$ be a commutative $k$-algebra. The $k$-Lie algebra $\Der_k(A)$ of $k$-linear derivations of $A$ is a Lie-Rinehart algebra with the identity map as anchor map.
We may define the ring of $k$-linear polynomial differential operators of $A$ denoted $\Diff_k(A)$. The ring $\Diff_k(A)$ has a canonical filtration $\Diff^i_k(A)$ give by the "order" of a differential operator. There
is an equality $\Diff^1_k(A) \cong A\oplus \Der_k(A)$, hence the Lie-Rinehart algebra $\Der_k(A)$ is a sub-$k$-Lie algebra and sub-$A$-module of $\Diff^1_k(A)$, where $\Diff^1_k(A)$ is the module of first order $k$-linear  differential operators of $A$. The algebra $\Diff^1_k(A)$ is an example of a D-Lie algebra, with the identity map and zero 2-cocycle defining it: There is an isomorphism

\[  \Diff^1_k(A) \cong \operatorname{D}_0^1(A) \]

of $k$-Lie algebras and $A\otimes_k A$-modules. Let $\alpha:=Id$ and $\pi: \Diff^1_k(A) \rightarrow \Der_k(A)$ be the projection map. Let $z:=Id_A \in \Diff^1_k(A)$ be the identity endomorphism of $A$.
There is a $k$-Lie bracket $[,]$ is defined on $\Diff^1_k(A)$ and one checks the 5-tuple $(\Diff^1_k(A), \alpha, \pi, [,], z)$ is a D-Lie algebra.

In the book \cite{beidar} the reader will find the notion of a \emph{differential Lie algebra} where $k\rightarrow A$ is an arbitrary unital map of commutative unital rings, and it seems
this notion is related to the notion of a pre-$\Dl$-Lie algebra. A differential Lie algebra is a generalization of a notion introduced in a paper of Baer from 1927. A differential $A/k$-Lie algebra
is a Lie-Rinehart algebra with a right $A$-module structure but it differs from the above definition of a D-Lie algebra.

I have not found the notion $\Dl$-Lie algebra in the litterature. It depends in a non-trivial way on the notion of a 2-cocycle for $\Der_k(A)$ and this notion was first introduced in Rineharts PhD-thesis \cite{rinehart}. Rineharts thesis was apparently the first paper to write down a proof of the fact that the differential $d^p_{\nabla}$ gave rise to a complex
\begin{align}
&\label{LRcomplex}  d^p_{\nabla}: \Hom_A(\wedge^p_A L ,(W,\nabla)) \rightarrow \Hom_A(\wedge^{p+1}_A L, (W,\nabla)). 
\end{align}
when $(W,\nabla)$ is a flat $L$-connection.
The paper \cite{rinehart} was the first paper to prove the PBW-theorem for $\U(A,L)$ when $L$ is a projective $A$-module and to systematically study the complex \ref{LRcomplex} and the cohomology 
$\H^i(L,(W,\nabla))$ and homology $\H_i(L,(W,\nabla))$ using $\U(A,L)$. For this reason many authors refer to the complex \ref{LRcomplex} as the \emph{Lie-Rinehart complex}. The paper
\cite{rinehart} is one of the most cited papers in this field.

\begin{lemma}\label{rightcomp} Let $(\tL, \ta, \tp, [,],\iD)$ be a $\Dl$-Lie algebra and let $u,v \in \tL$ and $c\in A$. It follows
\[ [u,vc]=c[u,v]+\tp(u)(c)v+\tp(u)\tp(v)(c)\iD.\]
\end{lemma}
\begin{proof} We get since $\iD$ is the central element in $\tL$ the following calculation:
\[ [u,vc]=[u,cv+\tp(v)(c)\iD]=c[u,v]+\tp(u)(c)v+\tp(u)\tp(v)(c)\iD .\]
The Lemma follows.
\end{proof}

\begin{example} Rings of differential operators and principal parts. \end{example}

Note: An $A\otimes_k A$-module $W$ that is finitely generated and projective as left and right $A$-module separately is called an \emph{$(A,A)$-vector bundle}. 
The module of principal parts $\P^l(E)$ is an $(A,A)$-vector bundle in many cases. There are examples where the left structure on $\P^l(E)$ is different from the right structure (see \cite{maa15}).
Similar results hold for the module of differential operators $\Diff^l(E,E)$. From \cite{maa15} we get the following example. Let $C:=\mathbb{P}^1$ be the projective line over a field of characteristic zero and let 
$\O(n)$ be the invertible sheaf with $n\in \mathbb{Z}$ an integer. The module of $l$'th order differential operators $\Diff^l(\O(n))$ from $\O(n)$ to $\O(n)$ has a left and right structure as $\O_C$-module and we get the following classification:

\begin{theorem} Let $k\geq1$ and $n\in \mathbb{Z}$ be integers. The following holds:
\[ \Diff^l(\O(n))^{right}\cong \O_{C}\oplus \O(l+1)^k \text{ for all $l\geq 1$ and $n\in \mathbb{Z}$.}\]
\[ \Diff^l(\O(n))^{left}\cong  \O(l)^{l+1} \text{ for all $1 \leq l\leq n$.}\]
\[ \Diff^l(\O(n))^{left}\cong \O_{C}^{n+1}\oplus \O(l+1)^{n-l} \text{ for all $n < l$ and $ l \geq 1$}\]
\end{theorem}
\begin{proof} The proof follows from \cite{maa15} since the sheaf of differential operators $\Diff^l(\O(n))$ is the dual of the sheaf of principal parts. If $1\leq k \leq n$ we get the following calculation:
\[ \Diff^l(\O(n))\cong \Hom_{\O_{C}}(\J^l(\O(n), \O(n)) \cong \Hom_{\O_{C}}( \O(n-l)^{l+1}, \O(n)) \cong \]
\[ \Hom_{\O_{C}}( \O(n-l ), \O(n))^{l+1} \cong   \Hom_{\O_{C}}( \O , \O(l))^{l+1} \cong  \O(l)^{l+1}.\]
The rest of the Theorem follows similarly.
\end{proof}

Hence for modules of differential operators $\Diff^l(\O(n))$ on the invertible sheaves $\O(n)$ on $C$, it follows the left and right $\O$-module structure differ in many cases. This phenomenon is
similar to the case for the module of principal parts, because $\Diff^l(\O(n))$ is defined as the dual of the module of principal parts. We may expect something similar to happen for a $\Dl$-Lie algebra
$\tL$, since it is a generalization of the first order module of differential operators $\Dl^1_0(A)$ on $A/k$.

For a line bundle $L\in \Pic(A)$ it follows the ring of differential operators $\Diff(L)$ is locally isomorphic to $\Diff(A)$ in the following sense: If $U:=\Spec(B)\subseteq \Spec(A)$ is an open subset where $L$ trivialize
it follows $\Diff(L)(U)\cong \Diff(B)$. Hence $\Spec(A)$ has an open affine cover $U_i:=\Spec(B_i)$ with $\Diff(L)(U_i)\cong \Diff(B_i)$. Here we abuse the notation and write 
$\Diff(A)$ when we speak of the sheaf on $\Spec(A)$ associated to $\Diff(A)$. It follows the module $\Diff^l(\O(n))$ is locally isomorphic to the module of differential operators $\Diff^l(\O_{C})$ of order $\leq l$.
One could try to generalize Definition \ref{dLie} to more general rings of differential operators on the form $\Diff(L)$ where $L$ is a line bundle or a projective $A$-module of rank $\geq 1$. 
One may replace
the $A\otimes_k A$-module $\D$ with $\Diff^1(L)$ where $L$ is a linebundle or a local twist of $\Diff^1(L)$. If $U_i:=\Spec(B_i)$ is an open cover where $L$ trivialize and $f_i\in \Z^2(\Der_k(B_i),B_i)$ is a 2-cocycle
we may locally define the $A\otimes_k A$-module $\Dl^1_{f_i}(B_i)$ and try to glue these to a global object $\Dl^1_f(L):=\Dl^1_{\{f_i\}}(L )$. Definition \ref{dLie} is functorial hence 
it makes sense to ask such a question.

\begin{example} Universal central extensions of Lie algebras and Lie-Rinehart algebras. \end{example}

Note: For Lie algebras and Lie-Rinehart algebras there is the notion of a \emph{universal central extension}. This is a new topic that is in development 
 (see the paper \cite{castiglioni} from 2014 for the case of Lie-Rinehart algebras).  One would like to give a construction of a "universal central extension" of a $\Dl$-Lie algebra. 

\begin{example} Connections in principal fiber bundles.\end{example}

The  constructions in this paper will in future work be generalized to give a theory valid for principal fiber bundles and 
connections in principal fiber bundles. This is work in progress (see  \cite{maa100}).

\begin{example} The $\Dl$-Lie algebra  associated to an Lie-Rinehart algebra. \end{example}

There is for every 2-cocycle a functorial way to construct a $\Dl$-Lie algebra from an Lie-Rinehart algebra.
Let $\alpha:L\rightarrow \Der_k(A)$ be a Lie-Rinehart algebra and let 
$f\in \Z^2(\Der_k(A),A)$ be a 2-cocycle. Let $\fal \in \Z^2(L,A)$ be the pull back 2-cocycle. Make the following definition:
\[ L(\fal):= Az\oplus L \]
with the following Lie bracket. Let $u:=az+x, v:=bz+y\in L(\fal)$. Define
\[ [u,v]:= (\alpha(x)(b)-\alpha(y)(a)+\fal(x,y),[x,y]) \in L(\fal).\]
Define the map
\[ \alpha_f:L(\fal)\rightarrow \D \]
by
\[ \alpha_f(az+x):=(a,\alpha(x))\in \D:=A\oplus \Der_k(A).\]
Define for $c\in A$
\[ cu:=c(az+x):=(ca)z+cx \]
and
\[ (az+x)c:=(ac+\alpha(x)(c))z+cx.\]
One checks that
\[ uc=cu+\alpha_f(u)(c) z \]
and $[z,u]=0$ hence the element $z\in L(\fal)$ is the canonical central element in $Z(L(\fal))$. One checks that the 5-tuple $(L(\fal), \alpha_f, \pi_f,[,], z)$ is a $\Dl$-Lie algebra in the sense
of Definition \ref{dLie}. A map of Lie-Rinehart algebra $\phi: (L,\alpha)  \rightarrow (L', \beta) $ gives in a canonical way a map of $\Dl$-Lie algebras
\[ \phi_f: L(\fal) \rightarrow L'(f^{\beta}) \]
defined by
\[ \phi_f(az+x):= az'+\phi(x).\]
Hence $\phi_f(z)=z'$. Hence we get a functor 
\[ F: \LR \rightarrow \underline{\Dl-\text{Lie}} \]
from the category $\LR$ of Lie-Rinehart algebras to the category of $\Dl$-Lie algebras  $\underline{\Dl-\text{Lie}}$   defined by
\[ F(L, \alpha):= (L(\fal), \alpha_f, \pi_f,[,], z).\]

\begin{example} $\Dl$-Lie algebras and connections.\end{example}

\begin{definition}
Let  $(\tL, \ta, \tp,[,],\iD)$ be a $\Dl$-Lie algebra with central element $\iD\in Z(\tL)$ and let $\psi \in \End_A(E)$ where $E$ is a left $A$-module.  An \emph{$(\tL, \psi)$-connection} is a left $A$-linear map
\[ \rho: \tL \rightarrow \End_k(E) \]
where the following holds for all $a\in A, u\in \tL$ and $e\in E$:
\[ \rho(u)(ae)=a\rho(u)(e)+\tp(u)(a)\psi(e).\]
Let $\Conn(\tL, \End)$ denote the category of all $(\tL,\psi)$-connection and morphisms of connections. Note: The endomorphism $\psi$ may vary. 

An \emph{$\tL$-connection} is an $A\otimes_k A$-linear map
$\rho: \tL\rightarrow \End_k(E)$ where $E$ is a left $A$-module. Let $\Mod(\tL)$ denote the category of all
$\tL$-connections and morphisms ot $\tL$-connections. Let $\Mod(\tL, Id)$ denote the category of all $\tL$-connections $(E,\rho)$ with $\rho(\iD)=Id_E$. It follows there is an inclusion of categories
\[ \Mod(\tL, Id) \subseteq \Mod(\tL).\]
\end{definition}

\begin{lemma} Let $\rho: \tL \rightarrow \End_k(E)$ be an $(\tL, \psi)$-connection. It follows for any element $x\in \tL$ the operator $\rho(x)$ is a differential operator of order one.
Hence if $f\in \Z^2(\Der_k(A),A)$ is a 2-cocycle and $\Dl^1_f(A)$ is the $\Dl$-Lie algebra associated to $f$, we get a canonical map
\[ \rho: \Dl^1_f(A) \rightarrow \Dl^1_0(A) \]
defined by
\[ \rho(a,x):=\phi_a + x \]
where $\phi_a:A\rightarrow A$ is multiplicaiton by $a$. The map $\rho$ is $A\otimes_k A$-linear and it follows the curvature $R_{\rho}$ satisfies
\[R_{\rho}(u,v)(c)=f(x,y)c ,\]
where $u:=(a,x), v:=(b,v)\in \Dl^1_f(A)$ and $c\in A$. Hence $\rho$ is a non-flat connection.
\end{lemma}
\begin{proof} The proof is an exercise.
\end{proof}.




\begin{lemma} Any $\tL$-connection $(E,\rho)$ is an $(\tL,\psi)$-connection for some $\psi\in \End_A(E)$, hence there is an inclusion of categories
\[ \Mod(\tL, Id) \subseteq \Mod(\tL) \subseteq \Conn(\tL, \End).\] 
\end{lemma}
\begin{proof} Let $\iD\in Z(\tL)$ be the canonical central element. It follows for any $v\in \tL$ and $c\in A$
\[ vc=cv +\tp(v)(c)\iD.\]
We get for any $\tL$-connection $\rho:\tL\rightarrow \End_k(E)$ the following:
\[ \rho(v)(ce)=\rho(vc)(e)=\rho(cv +\tp(v)(c)\iD)(e)=c\rho(v)(e)+\tp(v)(c)\rho(\iD)(e).\]
We get for any element $a\in A$
\[ \iD a =a\iD +\tp(\iD)(a)\iD =a\iD \]
since $\tp(\iD)=0$. It follows 
\[ \rho(\iD)(ae)=\rho(\iD a)(e)=\rho(a\iD)(e)=a\rho(\iD)(e)\]
hence $\rho(\iD):=\psi \in \End_A(E)$. It follows $\rho$ is an $(\tL, \psi)$-connection and the Lemma follows
\end{proof}

\begin{example} The universal ring of a $\Dl$-Lie algebra.\end{example}

\begin{definition} \label{universal} Given a $\Dl$-Lie algebra $(\tL, \ta,\tp,[,], \iD)$, define $\Uo(\tL)$ as in \cite{maa1}, Definition 3.3.
\end{definition}

Note: It follows from \cite{maa145} that $\Uo(\tL)$ is an associative unital ring which is functorial in the $\Dl$-Lie algebra $\tL$.

\begin{lemma} Let $f\in \Z^2(\Der_k(A),A)$ and let $(L,\alpha)$ be a Lie-Rinehart algebra. There is a covariant functor
\[  F_f: \LR \rightarrow \underline{Rings} \]
defined by
\[ F_f((L,\alpha)):= \Uo(L(\fal))    \]
where $L(\fal)$ is the $\Dl$-Lie algebra introduced in \cite{maa145}, Theorem 2.7.
Here $\underline{Rings}$ is the category of associative unital rings and maps of unital rings.
\end{lemma}
\begin{proof} Since $L(\fal)$ is a $\Dl$-Lie algebra the proof is clear from the above discussion. 
\end{proof}

\begin{theorem} \label{mainequivcat} Let $(\tL, \ta,\tp,[,],\iD)$ be a $\Dl$-Lie algebra. There is an exact equivalence of categories
\[ F: \Mod(\tL, Id)\cong \Mod(\Uo(\tL)) \]
preserving injective and projective objects.
\end{theorem}
\begin{proof} See \cite{maa145} for a proof.
\end{proof}

By Theorem \ref{mainequivcat} we may interpret the category $\Mod(\tL, Id)$ as a category of modules over an associative ring $\Uo(\tL)$. We may use Theorem \ref{mainequivcat}
to define Ext and Tor groups of arbitrary pairs of connections $(E,\rho),(E',\rho')$ in $\Mod(\tL, Id)$.

Note: In \cite{beidar} a \emph{differential $A/k$-Lie algebra} $L$  is introduced, and $L$ has a universal enveloping algebra $\U(L)$ with the following property: When the base ring $A$ is Noetherian and $L$ a finitely generated and projective $A$-module it follows $\U(L)$ is Noetherian. The ring $\U(L)$ is moreover almost commutative and has a PBW-basis when $L$ is projective as $A$-module.
The ring $\Uo(\tL)$ from Definition  \ref{universal} is not Noetherian in general and it is not almost commutative hence the notion $\Dl$-Lie algebra differs from the notion differential $A/k$-Lie algebra.

Let in the following  $(W,[,])$ be an $A$-Lie algebra where $W$ is an $A\otimes_k A$-module with $aw=wa$ for all $w\in W$ and $a\in A$. Assume
\begin{align}
&\label{exact1}0\rightarrow W \rightarrow \tilde{L_1} \rightarrow \tL \rightarrow 0 
\end{align}
is the exact sequence induced by a surjective map $p_1:\tilde{L_1}\rightarrow \tL$ of $\Dl$-Lie algebras $\tilde{L_1},\tL$ with $W:=ker(p_1)$. 
Let $(\tLone, \pi_1, \tilde{\pi}_1,[,], \iota_1)$ and $(\tL, \pi, \tilde{\pi}, [,], \iD)$ be the pairs of 5-tuples of elements corresponding to $\tLone$ and $\tL$.

Assume $\tL$ is projective as left $A$-module and that $s:\tL\rightarrow \tilde{L_1}$
is a left $A$-linear map with $p_1 \circ s =id_{\tL}$. Associated to the splitting $s$ we may defined the pair $(\nabla, \psi)$ as follows: 
Proposition \ref{LRalgebra} constructs the pair $(\nabla, \psi)\in S((\tL, \tp), (W,[,]))$ from the section $s$, where $(\tL, \tp)$ is the underlying Lie-Rinehart algebra of $\tL$. 
It follows $\nabla(x):=[s(x), -]$ and $\psi(x,y):=[s(x),s(y)]-s([x,y])$.
Let $(\tL(W, (\nabla,\psi)),\ta_W)$ be the Lie-Rinehart algebra constructed in Proposition \ref{LRalgebra}. Hence there is an equality of left $A$-modules $\tL(W,(\nabla,\psi))\cong W\oplus \tL$.

Let $\iD_s:=(\iD_1-s(\iD), \iD)\in W\oplus \tL$ where $\iD_1\in Z(\tL_1)$ and $\iD\in Z(\tL)$ are the central elements of $\tL_1$ and $\tL$.

\begin{lemma} The element $\iD_s$ is a central element in the $k$-Lie algebra $\tL(W,(\nabla,\psi))$.
\end{lemma}
\begin{proof} Let $u:=(w,x)\in W\oplus \tL:=\tL(W,(\nabla,\psi))$. We get
\[ [u,\iD_s]:=[(w,x),(\iD_1-s(\iD),\iD)]:= \]
\[ ([w,\iD_1-s(\iD)]+[s(x), \iD_1-s(\iD)]-[s(\iD),x]+[s(x), s(\iD)]-s([x,\iD]), [x,\iD])=\]
\[ ([w,\iD_1]-[w,s(\iD)]+[s(x), \iD_1]-[s(x), s(\iD)]-[s(\iD),w]+[s(x), s(\iD)]+ \]
\[  -s([x,\iD]), [x,\iD])=\]
\[ ( -[w,s(\iD)] -[s(\iD),w] -[s(x),s(\iD)]+[s(x), s(\iD)], 0) =0\]
and it follows $\iD_s$ is a central element in $\tL(W,(\nabla, \psi))$.
\end{proof}

\begin{lemma}  \label{centralexistence} Let $(\tL, \ta, \tp,[,],\iD)$ be a $\Dl$-Lie algebra and let $(W,[,])$ be an $A\otimes_k A$-module and $A$-Lie algebra with $aw=wa$ for all $a\in A, w\in W$.
Let $ \iD_s:=(w_1,x_1)\in W\oplus \tL$. Let $(\nabla, \psi) \in S(\tL, (W,[,]))$ and let $(\tL(W,(\nabla,\psi)), \tp_W)$ be the Lie-Rinehart algebra introduced in Proposition \ref{LRalgebra}.
Let $p:\tL(W,(\nabla,\psi)) \rightarrow \tL$ be the canonical projection map.
It follows the element $\iD_s$  is a central element in $\tL(W,(\nabla, \psi))$ with $p(\iD_s)=\iD$ if and only if $x_1=\iD$ and $w_1$ satisfies the following equations:
\begin{align}
&\label{central1} \nabla(\iD)=ad(-w_1) \\
&\label{central2} \psi(\iD,-)= d^0(w_1)
\end{align}
\end{lemma}
\begin{proof} If $\iD_s:=(w_1,x_1)\in W\oplus \tL$ it follows $p(\iD_s)=\iD$ if and only if $x_1=\iD$ hence $\iD_s=(w_1,\iD)$. 
One checks that $\iD_s$ is a central element if and only if equations \ref{central1} and \ref{central2} holds, and the 
Lemma is proved.
\end{proof}

Let $(\tL, \ta, \tp,[,],\iD)$ be a $\Dl$-Lie algebra and let $(W,[,])$ be an $A\otimes_k A$-module and $A$-Lie algebra with $aw=wa$ for all $a\in A, w\in W$.
Let $(\nabla, \psi)\in S(\tL, (W,[,]))$. Let $\iD_s:=(w_1,w_2)\in W\oplus \tL:=\tL(W,(\nabla,\psi))$. There are canonical maps
\[ \ta_W: \tL(W,(\nabla, \psi)) \rightarrow \Dl^1_f(A) \]
defined by
\[ \ta_W(w,x):=\ta(x)\]
and
\[ \tp_W:\tL(W,(\nabla,\psi)) \rightarrow \Der_k(A) \]
defined by
\[ \tp_W(w,x):=\tp(x).\]
By Proposition \ref{LRalgebra} it follows the pair $(\tL(W,(\nabla, \psi)), \tp_W)$ is an Lie-Rinehart algebra and the canonical sequence
\begin{align}
&\label{LRext} 0 \rightarrow W\rightarrow \tL(W,(\nabla,\psi)) \rightarrow \tL \rightarrow 0
\end{align}
is an extension of Lie-Rinehart algebras.

\begin{theorem} \label{mainDLIE} Let $(\tL,\ta, \tp,[,],\iD)$ be a $\Dl$-Lie algebra and let $(W,[,])$ be an $A\otimes_k A$-module and $A$-Lie algebra with $aw=wa$ for all $a\in A, w\in W$. Let $(\nabla, \psi)\in
S(\tL,(W,[,]))$. Let $\iD_s:=(w_1,x_1)\in \tL(W(\nabla, \psi))$. There is a structure as $\Dl$-Lie algebra on 
$W\oplus \tL$ making the 5-tuple $(W\oplus \tL, \ta_W, \tp_W,[,],\iD_s)$ into a $\Dl$-Lie algebra with \ref{LRext} an extension of $\Dl$-Lie algebras if and only if the element $\iD_s$ 
satisfies $x_1:=\iD$ and $w_1$ satisfies equations \ref{central1} and \ref{central2}.
\end{theorem}
\begin{proof} Let $\iD_s:=(w_1,x_1)\in \tL(W,(\nabla,\psi))$. By Lemma \ref{centralexistence} it follows $\iD_s$ is a central element with $p(\iD_s)=\iD$ if and only if $x_1=\iD$ and $w_1$
satisfies equations \ref{central1} and \ref{central2}. Assume $\iD_s=(w_1,\iD)$ is a central element. Let $c\in A$ and define the following right $A$-module structure on the element $u:=(w,x)\in W\oplus \tL$:
\begin{align}
&\label{rightA} uc:=cu+\tp_W(u)(c)\iD_s
\end{align}
where
\[ \tp_W: \tL(W,(\nabla, \psi)) \rightarrow \Der_k(A) \]
is the canonical map defined by $\tp_W(w,x):= \tp(x).$ The claim is that  Definition \ref{rightA} defines a right $A$-module structure on $W\oplus \tL$ with $(c'u)c=c'(uc)$ for all $c,c'\in A$.
Let us verify the axioms of an $A\otimes_k A$-module structure on $\tL(W,(\nabla, \psi))$. One easily checks that for any elements $u,v\in W\oplus \tL$ we get
\[ (u+v)c=uc+vc \]
and
\[ u 1 =u \]
where $ 1\in A$ is the multiplicative unit.
We also get 
\[ (cu)c'=c(uc').\]
Let us verify that $(uc)c'=u(cc')$. We get
\[ u(cc')=(cc')u+\tp_W(u)(cc')\iD_s.\]
One checks that $\iD_sc'=c'\iD_s$.
We get
\[ (uc)c'=(cu+\tp_W(u)(c)\iD_s)c = (cu)c'+\tp_W(u)(c)\iD_s)c'=\]
\[ c'(cu)+\tp_W(cu)(c')\iD_s+ c'(\tp_W(u)(c)\iD_s)+\tp_W(\tp_W(u)(c))(c')\iD_s)(c')\iD_s =\]
\[ (c'c)u+c\tp_W(u)(c')\iD_s+ c'\tp_W(u)(c)\iD_s =  \]
\[ (cc')u+\tp_W(u)(cc')\iD_s=u(cc')\]
since $\tp_W(u)\in \Der_k(A)$. Hence $u(cc')=(uc)c'$. 
One checks that the 5-tuple $(W,\oplus \tL, \ta_W, \tp_W,[,],\iD_s)$ is a $\Dl$-Lie algebra and that \ref{LRext} is an exact sequence of $\Dl$-Lie algebras and the Theorem follows.
\end{proof}


\begin{definition} \label{conditionB} Let $(\tL, \ta, \tp,[,],D)$ be a $\Dl$-Lie algebra and let $(W,[,])$ be an $A\otimes_k A$-module and $A$-Lie algebra with $aw=wa$ for all $a\in A, w\in W$.
Let $D(\tL, (W,[,]))$ be the set of triples $(\nabla, \psi, w_1) $ with $(\nabla, \psi)\in S(\tL,(W,[,]))$ and where $w_1\in W$ is an element
satisfying equations \ref{central1} and \ref{central2}. A triple $(\nabla, \psi, w_1)\in D(\tL, (W,[,]))$ is said to satisfy \emph{condition B}.
\end{definition}

\begin{definition} \label{equivalenceD} Let $(\tL, \ta, \tp,[,],D)$ be a $\Dl$-Lie algebra and let $(W,[,])$ be an $A\otimes_k A$-module and $A$-Lie algebra with $aw=wa$ for all $a\in A, w\in W$.
Define the following relation on the set $D(\tL, (W,[,]))$: Given two triples $s:=(\nabla,\psi, \tilde{w}), s':=(\nabla', \psi',\tilde{w}')$ in $D(\tL, (W,[,]))$. We say $s\equiv s'$ if and only if there is a map
$b\in \Hom_A(\tL,W)$ with the property that
\begin{align}
&\label{ec1} \nabla'(x)=\nabla(x)+[b(x),-]  \\
&\label{ec2} \psi'(x,y)=\psi(x,y)+d^1_{\nabla}(b)(x,y)+[b(x),b(y)] \\
 &\label{ec3}\tilde{w}'=\tilde{w}-b(\iD) 
\end{align}
where $\iD\in Z(\tL)$ is the canonical central element.
\end{definition}

\begin{lemma} \label{equivS} The relation $\equiv$ from Definition \ref{equivalenceD} is an equivalence relation on $D(\tL,(W,[,]))$.
\end{lemma}
\begin{proof} Let $b=0$ It follows
\[ \nabla(x)=\nabla(x) +[b(x),-] \]
\[ \psi(x,y)=\psi(x,y)+d^1_{\nabla}(b)(x,y)+[b(x),b(y)] \]
and
\[ \tilde{w}=\tilde{w}-b(\iD) \]
hence if $s:=(\nabla, \psi, \tilde{w})$ it follows $s \equiv s$. 

Assume $s':=(\nabla',\psi',\tilde{w}')$ is another element with $s\equiv s'$ via an element $b\in \Hom_A(\tL,W)$. We get the following calculation:
\[ \nabla(x)=\nabla'(x)-[b(x),-] \]
hence if we let $a:=-b$ we get
\[\nabla(x)=\nabla'(x)+[a(x),-] .\]
We get
\[ \psi(x,y)+d^1_{\nabla}(b)(x,y)+[b(x),b(y)]=\]
\[ \psi(x,y)+\nabla(x)(b(y))-\nabla(y)(b(x))-b([x,y])+[b(x),b(y)] =\]
\[ \psi(x,y)+\nabla'(x)(b(y)) +[a(x),b(y)]-\nabla'(y)(b(x))-[a(y),b(y)]-b([x,y])+[b(x),b(y)]=\]
\[ \psi(x,y)-(\nabla'(x)(a(y))-\nabla'(y)(a(x))a([x,y])) -[b(x),b(y)]=\]
\[ \psi(x,y)-d^1_{\nabla'}(a)(x,y)-[a(x),a(y)].\]
Hence we get
\[ \psi'(x,y)=\psi(x,y)-d^1_{\nabla'}(a)(x,y)-[a(x),a(y)]\]
and it follows 
\[ \psi(x,y)=\psi'(x,y)+d^1_{\nabla'}(a)(x,y)+[a(x),a(y)].\]
Since $\tilde{w}'=\tilde{w}-b(\iD)$ it follows
\[ \tilde{w}=\tilde{w}'+b(\iD)=\tilde{w}'-a(\iD) \]
and we have shown that $s' \equiv s$ and the relation is symmetric.
 Assume $s'':=(\nabla'', \psi'', \tilde{w}'')$ and $s\equiv s'$ via an element $b\in \Hom_A(\tL,W)$. Moreover $s' \equiv s''$ via an element $a\in \Hom_A(\tL,W)$. It follows from Lemma \ref{relation}
there is an equality
\[ \nabla''(x)=\nabla(x)+[(a+b)(x),-] \]
and
\[ \psi''(x,y)=\psi(x,y)+d^1_{\nabla}(a+b)(x,y)+[(a+b)(x),(a+b)(y)].\]
We get
\[ \tilde{w}''=\tilde{w}'-a(\iD)=\tilde{w}-b(\iD)-a(\iD)=\tilde{w}-(a+b)(\iD) \]
hence $s\equiv s''$. It follows the relation $\equiv$ is an equivalence relation on $D(\tL, (W,[,]))$ and the Lemma is proved.
\end{proof}

\begin{definition} Let $(\tL, \ta,\tp,[,],\iD)$ be a $\tL$-Lie algebra and $(W,[,])$ an $A$-Lie algebra with $aw=wa$ for all $a\in A$ and $w\in W$.
Define the quotient set
\[   \ext^1(\tL, (W,[,])):= D(\tL,(W,[,]))/\equiv \]
where $\equiv$ is the equivalence relation from definition  \ref{equivalenceD}.
\end{definition}

Note: Every element in $D(\tL,(W,[,]))$ may be represented by a triple $(\nabla,\psi, w)$ satisfying condition B.


\begin{definition} Let $(\tL, \ta, \tp,[,],\iD)$ be a $\Dl$-Lie algebra and let $(W,[,])$ be an $A\otimes_k A$-module and $A$-Lie algebra with $aw=wa$ for all $a\in A, w\in W$.
Let $(\nabla, \psi, \tilde{w})\in D(\tL, (W,[,]))$.
Let the 5-tuple $(\tL(W, (\nabla, \psi, \tilde{w})), \ta_W, \tp_W,[,],\iD_s)$ denote the $\Dl$-Lie algebra constructed in Theorem \ref{mainDLIE}.
\end{definition}

There is by Theorem \ref{mainDLIE} for every triple $(\nabla, \psi, \tilde{w}) \in D(\tL, (W,[,]))$ an extension
\[ 0 \rightarrow W \rightarrow  \tL(W,(\nabla, \psi,\tilde{w})) \rightarrow \tL \rightarrow 0 \]
of $\Dl$-Lie algebras.

\begin{definition}
Fix an $A$-Lie algebra and $A\otimes_k A$-module $(W,[,])$ with $aw=wa$ for all $a\in A$ and $w\in W$, and a $\Dl$-Lie algebra $(\tL, \ta,\tp,[,], \iD)$.
Let $T(\tL, (W,[,]))$ denote the set of exact sequences of $\Dl$-Lie algebras
\[ 0 \rightarrow W \rightarrow \tilde{L}_1 \rightarrow \tL \rightarrow 0 .\]
\end{definition}

Define two extensions $p_1: \tilde{L}_1 \rightarrow \tL$ and $p_2:\tilde{L}_2 \rightarrow \tL$ to be equivalent if and only if there is an isomorphism $\rho: \tilde{L}_1 \rightarrow \tilde{L}_2$ of $\Dl$-Lie algebras
making all diagrams commute. Write $(\tilde{L}_1,p_1) \equiv (\tilde{L}_2, p_2)$. Since $W=ker(p_1)=ker(p_2)$ we drop $W$ from the notation.
It follows $ \equiv$ is an equivalence relation on $T(\tL, (W,[,]))$ since $\rho$ is an isomorphism. 

\begin{definition}
Let $\Ext^1(\tL, (W,[,]))$ denote the quotient set $T(\tL, (W,[,]))/\equiv$.
\end{definition}

\begin{theorem} \label{maintheorem} Let $(\tL, \tp, \tp, [,], \iD)$ be a $D$-Lie algebra where $\tL$ is projective as left $A$-module and let $(W,[,])$ be an $A$-Lie algebra with $aw=wa$ for all $a\in A$ and $w\in W$.
Given a surjective map $p:\tilde{L}'\rightarrow \tL$ of $\Dl$-Lie algebras with left $A$-linear section $s$ and define the connection $\nabla_s: \tL \rightarrow \Der_k(W)$ by $\nabla_s(x):=[s(x),-]$ and 2-cocycle
$\psi_s: \wedge^2 \tL \rightarrow W$ by $\psi_s(x,y):=[s(x),s(y)]-s([x,y])$. Define the element $w_s:=\iD'-s(\iD)\in W$ where $\iD'\in Z(\tilde{L}')$ is the central element.
Define the following map
\[ F: \Ext^1(\tL, (W,[,])) \rightarrow \ext^1(\tL, (W,[,])) \]
by
\[ F(\tilde{L}',p)  :=(\nabla_s, \psi_s, w_s)\operatorname{mod} \equiv.\]
The map $F$ is a well defined one-to-one correspondence of sets.
\end{theorem}
\begin{proof} Assume $p:\tilde{L}'\rightarrow \tL$ is a surjective map of $\Dl$-Lie algebras with splitting as left $A$-modules $s$. Let $F(\tilde{L}', p):=(\nabla_s, \psi_s, w_s)$. Assume
$s':=s+b$ is another splitting with $b\in \Hom_A(\tL,W)$. It follows
\[ \nabla_{s'}(x)=\nabla_s(x)+[b(x), -] \]
\[ \psi_{s'}(x,y)=\psi_s(x,y)+d^1_{\nabla}(x,y)+[b(x),b(y)] \]
and
\[ w_{s'}=\iD_s -s'(\iD):=\iD_s-s(\iD)-b(\iD)=w_s-b(\iD) \]
hence
\[  (\nabla_s, \psi_s, w_s)\equiv (\nabla_{s'}, \psi_{s'}, w_{s'}) \]
and $F$ is a well defined map of sets. Assume $p_1:\tilde{L_1}\rightarrow \tL$ and $p_2:\tilde{L_2}\rightarrow \tL$ are two surjective maps of $\Dl$-Lie algebras with 
\[ F(\tilde{L_1}, p_1)=F(\tilde{L_2}, p_2)\operatorname{mod}(\equiv) .\]
It follows there is an element $b\in \Hom_A(\tL, W)$ with 
\[ \nabla_{s'}(x)=\nabla_s(x)+[b(x), -] \]
\[ \psi_{s'}(x,y)=\psi_s(x,y)+d^1_{\nabla}(x,y)+[b(x),b(y)] \]
and
\[ w_{s'}=w_s-b(\iD) .\]
It follows there is an isomorphism
\[ \rho: \tL(W, (\nabla_{s'},\psi_{s'}, w_{s'})) \cong \tL(W, (\nabla_{s},\psi_{s}, w_{s}))  \]
of $\Dl$-Lie algebras inducing an isomorphism $\tilde{\rho}: \tilde{L_1}\rightarrow \tilde{L_2} $ of extensions of $\tL$. It follows $F$ is an injective map.
Let $s:=(\nabla,\psi, w)\in D(\tL, (W,[,]))$ be an element and define the $\Dl$-Lie algebra $(\tL(W, (\nabla, \psi,w)), \ta_s, \tp_s,[,], \iD_s)$ where  $\iD_s:=(w,\iD)$ as in Theorem \ref{mainDLIE}.
The canonical map
\[ p: \tL(W, (\nabla_s, \psi_s, w_s)) \rightarrow \tL \]
is a surjective map of left $A$-modules and since $\tL$ is projective it follows there is a left  $A$-linear  section $t$ of $p$ with $t(x):=(b(x),x)$ where $b\in \Hom_A(\tL, W)$. We
get
\[ \nabla_t(x)=\nabla(x)+[b(x), -] \]
\[ \psi_t(x,y)=\psi(x,y)+d^1_{\nabla}(b)(x,y)+[b(x),b(y)] \]
and
\[ w_t=w-b(\iD) \]
hence 
\[  (\nabla, \psi, w)\equiv (\nabla_t, \psi_t, w_t) \]
and the map $F$ is surjective. The Theorem is proved.
\end{proof}

\begin{example} The $\Dl$-Lie algebra associated to a connection.\end{example}

Let $(L,\alpha)$ be an Lie-Rinehart algebra and let $\rho:L\rightarrow \End_k(E)$ be a connection. 
Assume $E$ is a finitely generated and projective $A$-module of rank $r$ where the characteristic of $k$ does not divide $r$. 

There is a canonical isomorphism of left $A$-modules

\[ f: E^*\otimes_A E \rightarrow \End_A(E) \]

defined by
\[ f(\phi \otimes e)(x):=\phi(x)e.\]
There is an $A$-linear trace map
\[ tr: E^*\otimes_A E \rightarrow A \]
defined by
\[ tr(\phi \otimes e):= \phi(e). \]
Let $\tilde{\rho}:L\rightarrow \Der_k(\End_A(E))$ be defined by $\tilde{\rho}(x)(\phi):=[\rho(x),\phi]$. The connection $\tilde{\rho}$ is the tensor product of $\rho^*$ and $\rho$.

\begin{lemma} \label{dtr} Let $\End_A(E)^{tr}$ be the $A$-module of endomorphisms $\phi$ with
$tr(\phi)=0$.  Let $\End^{tr}_A(E):=\End_A(E)/\End_A(E)^{tr}$. The following holds: For any element $\phi \in \End_A(E)$ and $x \in L$ we get
\begin{align}
&\label{end1} tr(\tilde{\rho}(x)(\phi))=\alpha(x)(tr(\phi)) \\
&\label{end2}\text{If $tr(\phi)=0$ it follows $tr(\tilde{\rho}(x)(\phi))=0$.}\\
&\label{end3}tr([\phi,\psi])=0\\
&\label{end4}\text{$\End_A(E)^{tr}$ is an $A$-module and $A$-Lie algebra}\\
&\label{end5}\text{There is a connection $\rho^{tr}:L\rightarrow \Der_k(\End^{tr}_A(E))$}\\
&\label{end6}R_{\rho^{tr}}(x,y)(\overline{\phi})=[R_{\rho}(x,y),\overline{\phi}]
\end{align}
\end{lemma}
\begin{proof}
For an element $\phi\in \Hom_A(E,A)$ we defined the dual connection $\nabla^*$ as follows:
\[ \rho^*(x)(\phi):= \alpha(x)\circ \phi -\phi \circ \rho(x).\]
we get if $\phi:= \sum_i \phi_i \otimes e_i \in E^* \otimes_A E$ the following calculation:
\[ tr( \tilde{\rho}(x)(\phi))= tr(\sum_i \rho^*(\phi_i)\otimes e_i +\phi_i \otimes \rho(x)(e_i)) =\]
\[ tr(\sum_i (\alpha(x)\circ \phi_i -\phi_i\circ \rho(x) )\otimes e_i +\phi_i \otimes \rho(x)(e_i)) =\]
\[ \sum_i \alpha(x)(\phi_i(e_i))-\phi_i(\rho(x)(e_i))+\phi_i(\rho(x)(e_i))=  \]
\[ \sum_i \alpha(x)(\phi_i(e_i))=\alpha(x)(\sum_i \phi_i(e_i))=\alpha(x)(tr(\phi)) .\]
It follows \ref{end1}  is proved. If $tr(\phi)=0$ it follows $tr(\tilde{\rho}(x)(\phi))=\alpha(x)(tr(\phi))=x(0)=0$ hence \ref{end2} is proved. If $\phi:= \sum \phi_i \otimes e_i$ and $\psi:=\sum \psi_j\otimes f_j$ in 
$E^*\otimes_A E$ we get
\[\phi \circ \psi := \sum \phi_i \otimes \psi_j(e_i)f_j \]
hence
\[ tr(\phi \circ \psi)=\sum \phi_i(f_j)\psi_j(e_i).\]
Similarly we get
\[ tr(\psi \circ \phi)=tr(\sum \psi_j \otimes \phi_i(f_j)e_i)=\sum \psi_j(e_i)\phi_i(f_j) \]
hence
\[ tr(\phi \circ \psi)=tr(\psi \circ \phi).\]
It follows $tr([\phi,\psi])=0$ and \ref{end3} is proved. Claim \ref{end4} follows from \ref{end3}. We get a well defined connection
\[ \rho^{tr}:L\rightarrow \Der_k(\End^{tr}_A(E)) \]
with 
\[ R_{\rho^{tr}}(x,y)(\overline{\phi})=[R_{\rho}(x,y),\overline{\phi}] \]
and the Lemma follows.
\end{proof}

From Lemma \ref{dtr}, Claim \ref{end5} and \ref{end6} we may define the Lie-Rinehart algebra $\End(L,E):=\End_A^{tr}(E)\oplus L$ with the following Lie product: Let $u:=(\phi,x), v:=(\psi,y)\in \End(L,E)$ and define
\[ [u,v]:=([\phi,\psi]+\tilde{\rho}(x)(\psi)-\tilde{\rho}(y)(\phi)+R_{\rho}(x,y),[x,y]). \]
Here 
\[ R_{\rho}(x,y):=[\rho(x),\rho(y)]-\rho([x,y]) \in \End_A(E).\]
The Lie product is well defined on $\End_A^{tr}(E)$ from Lemma \ref{dtr}, claim \ref{end4}.
Define the map
\[ \beta: \End(L,E)\rightarrow \Der_k(A) \]
by
\[ \beta(\phi,x):=\alpha(x).\]
It follows the pair
\[ (\End(L,E), \beta) \]
is an Lie-Rinehart algebra with the following left $A$-module structure:
Let $c\in A$ and define the following left $A$-module structure on $\End(L,E)$:
\[ cu:=c(\phi,x):=(c\phi,cx) \]
Define the following right $A$-module struture:
\[ c:=(\phi,x)c:=(\phi c +\alpha(x)(c)I, cx) \]
where $I:=Id_E \in \End_A(E)$. One checks the multiplication is well defined.
Let $\iD:=(I,0)\in \End(L,E)$. It follows $\iD\in Z(\End(L,E))$ is a central element with the following property:
\[ uc=(\phi,x)c=(\phi c +\alpha(x)(c)I,cx)=c(\phi,x)+\alpha(x)(c)(I,0)=cu+\alpha(x)(c)\iD.\]
It follows
\[ uc=cu+\beta(u)(c)\iD\]
for all $u\in \End(L,E)$ and $c\in A$.
Since $\End(L,E)$ is a Lie-Rinehart algebra it follows 
\[ [u,cv]=c[u,v]+\beta(u)(c)v \]
for all $u,v\in \End(L,E)$ and $c\in A$.  Hence $\End(L,E)$ is an $A\otimes_k A$-module, a $k$-Lie algebra and the map $\beta$ is a map of $A\otimes_k A$-modules and $k$-Lie algebras.
It follows the 4-tuple $(\End(L,E), \beta,[,],\iD)$ is a pre-$\Dl$-Lie algebra.

The curvature $R_{\rho}\in \Z^2(L,\End_A(E))$  and since $E$ is projective of rank $r$ it follows the trace $f:=tr(R_{\rho})\in \Z^2(L,A)$. Hence we may define
the $k$-Lie algebra $L(f)$ associated to the 2-cocycle $f:=tr(R_{\rho})$. 

Note: We define the central element $z$ in $L(f)$ as follows: $z:=(r,0)=r(1,0)$ in $L(f)$.

We define the following left and right $A$-module structure on $L(f)$. Let $u:=(a,x)\in L(f)$ and let $c\in A$:
\[ cu=c(a,x):=(ca,cx) \]
and
\[uc=(a,x)c:=(ac+\alpha(x)(a)r,cx).\]
It follows
\[ uc=(ac+\alpha(x)(c)r,cx)=c(a,x)+\alpha(x)(c)(r,0) \]
hence
\[  uc= cu+\alpha_f(u)(c)z \]
where $\alpha_f: L(f) \rightarrow \Der_k(A) $ is the canonical projection map. It follows the 4-tuple
\[ (L(f), \alpha_f,[,],z) \]
is a pre-$\Dl$-Lie algebra. define the map
\[ \gamma: \End(L,E)\rightarrow L(f)\]
by
\[ \gamma(\phi,x):=(tr(\phi),x) .\]

Note: The $A$-Lie algebra $\End_A^{tr}(E)$ satisfies $a\phi=\phi a$ for all $\phi\in \End_A^{tr}(E)$ and $a\in A$.

\begin{theorem} \label{mainext} The 4-tuple $(\End(L,E),  \beta,[,],\iD)$ is a pre-$\Dl$-Lie algebra.  If $k$ contains a field of characteristic zero it follows 
the sequence
\[ 0 \rightarrow \End_A^{tr}(E) \rightarrow \End(L,E) \rightarrow L(f) \rightarrow 0\]
is a non-abelian extension of pre-$\Dl$-Lie algebras of $L(f)$ with the $A$-Lie algebra $(\End_A^{tr}(E),[,])$. 
\end{theorem}
\begin{proof} The left $A$-module $\End(L,E)$ is an Lie-Rinehart algebra by section one of this paper. We define $\iD:=(I,0)$ and it follows $\iD$ is a central element.
we see that $\gamma(I,0)=(tr(I),0)=(r,0)=z\in \Dl$ which is the central element in $\Dl$. It follows $\beta(\iD)=0$ in $\Der_k(A)$.
We get
\[ uc=cu+\beta(u)(c)\iD \]
for all $u\in \End(L,E)$ and $c\in A$. Since $\End(L,E)$ is a Lie-Rinehart algebra it follows 
\[ [u,cv]=c[u,v]+\beta(u)(c)v \]
for all $u,v\in\End(L,E)$ and $c\in A$. The map
\[ \gamma:\End(L,E) \rightarrow L(f) \]
is by definition a map of $A\otimes_k A$-modules. Assume $u:=(\phi,x), v:=(\psi,y)\in \End(L,E)$. We get from Lemma \ref{dtr}
\[ \gamma([u,v])=\gamma([\phi,\psi]+\tilde{\rho}(x)(\psi)-\tilde{\rho}(y)(\phi)+R_{\rho}(x,y),[x,y])=\]
\[ (tr([\phi,\psi])+tr(\tilde{\rho}(x)(\psi))-tr(\tilde{\rho}(y)(\phi))+tr(R_{\rho}(x,y)),[x,y]) =\]
\[ ( \beta(x)(tr(\psi))-\beta(y)(tr(\phi))+tr(R_{\rho}(x,y)),[x,y])=\]
\[ [(tr(\phi),x),( tr(\psi),y)] =[\gamma(u),\gamma(v)] \]
hence the map $\gamma$ is a map of $k$-Lie algebras. It follows $\gamma$ and $\beta$ are maps of $A\otimes_k A$-modules and $k$-Lie algebras.
If $k$ contains a field of characteristic zero it follows the canonical sequence
\[ 0\rightarrow \End_A^{tr}(E) \rightarrow \End(L,E) \rightarrow L(f) \rightarrow 0\]
is an exact sequence of left and right $A$-modules and the Proposition follows.
\end{proof}

Hence extensions of pre-$\Dl$-Lie algebras and $\Dl$-Lie algebras arise in a natural way when we study classical connections as Theorem \ref{mainext} shows.

\begin{example}  The non-abelian extension associated to a non-flat connection on a $\Dl$-Lie algebra.\end{example}

Let $k$ be the complex numbers (or a field of characteristic zero) and let $R:=k[x_i,1/x_i], (E,\rho)$ be the ring and flat connection from Example \ref{exnonab}.
Let $L_R:=\Der_k(R)$ and consider an element $\overline{f}\in \H^2(L_R,R)$. By Lemma \ref{derham} we may write
\[ \overline{f}=\sum_{1\leq i <j \leq n}\frac{\alpha_{ij}}{x_ix_j}dx_i\wedge dx_j \]
with $\alpha_{ij}\in k.$ We may lift $\overline{f}$ to an element $f\in \Z^2(L_R,R)\subseteq \Hom_R(\wedge^2 L_R,R)$. 

\begin{lemma} \label{2cocycle} Let $x:= \sum_{i}a_{i}\partial_{x_i}$ and $y:=\sum_{i}b_{i}\partial_{x_i}$ with $x,y\in  L_R$.
It follows
\[ f(x \wedge y)= \sum_{i<j} (a_ib_j-a_jb_i)\frac{\alpha_{ij}}{x_ix_j} \in R.\]
\end{lemma}
\begin{proof} The proof is an exercise.
\end{proof}

We may construct the $\Dl$-Lie algebra $L_R(f)$ associated to to 2-cocycle $f$.
We get from $(E,\rho)$ a connection
\[ \tilde{\rho}: L_R(f) \rightarrow\End_k(E) \]
defined by
\[ \tilde{\rho}(az+x):= aId_E+\rho(x).\]
Let $u:=(a,x), v:=(b,y)\in L_R(f)$. Its curvature is given by the formula
\[ R_{\tilde{\rho}}(u,v)=f(x\wedge y)Id_E= \sum_{i<j}(a_ib_j-a_jb_i)\frac{\alpha_{ij}}{x_ix_j} Id_E \in \End_A(E).\]
Hence the connection $\tilde{\rho}$ is non-flat. We get a non-abelian extension of $\Dl$-Lie algebras
\begin{align}
&\label{extnonabelian} 0 \rightarrow \End_A(E) \rightarrow \End(L(f),(E,\tilde{\rho})) \rightarrow L(f) \rightarrow 0.
\end{align}
Since $\tilde{\rho}$ is non-flat it follows the sequence \ref{extnonabelian} is non-split.

\textbf{Acknowledgements.} Thanks to Pierre Deligne and Marat Rovinsky for comments on an earlier version of this paper and families of connections. 
Thanks to Ernest Vinberg for the reference \cite{beidar}.

\end{document}